\documentclass[reqno]{amsart}

\usepackage{amsthm, amsmath, amssymb, amsfonts, graphicx}

\usepackage[english]{babel}

\usepackage{bbm}
\usepackage[foot]{amsaddr}

\usepackage[inline]{enumitem}

\usepackage{float}
\restylefloat{table}

\usepackage{graphicx}
\usepackage{caption}
\usepackage{subcaption}

\usepackage[colorlinks=true, pdfstartview=FitV, linkcolor=blue,
            citecolor=blue, urlcolor=blue]{hyperref}
\usepackage[usenames]{color} 

\usepackage[round]{natbib}

\newtheorem{theorem}{Theorem}
\newtheorem{proposition}[theorem]{Proposition}
\newtheorem{lemma}[theorem]{Lemma}
\newtheorem{corollary}[theorem]{Corollary}

\theoremstyle{definition}

\theoremstyle{remark}
\newtheorem{remark}[theorem]{Remark}

\numberwithin{equation}{section}
\numberwithin{theorem}{section}

\DeclareMathOperator*{\argmin}{arg\,min} 

\DeclareFontFamily{U}{matha}{\hyphenchar\font45}
\DeclareFontShape{U}{matha}{m}{n}{
      <5> <6> <7> <8> <9> <10> gen * matha
      <10.95> matha10 <12> <14.4> <17.28> <20.74> <24.88> matha12
      }{}
\DeclareSymbolFont{matha}{U}{matha}{m}{n}
\DeclareFontSubstitution{U}{matha}{m}{n}

\DeclareFontFamily{U}{mathx}{\hyphenchar\font45}
\DeclareFontShape{U}{mathx}{m}{n}{
      <5> <6> <7> <8> <9> <10>
      <10.95> <12> <14.4> <17.28> <20.74> <24.88>
      mathx10
      }{}
\DeclareSymbolFont{mathx}{U}{mathx}{m}{n}
\DeclareFontSubstitution{U}{mathx}{m}{n}

\DeclareMathDelimiter{\vvvert}{0}{matha}{"7E}{mathx}{"17}

\title{Long-run impulse control with generalised discounting}
\author{Damian Jelito$^{\ast,\dagger}$}

\email{djelito@impan.pl}

\author{\L{}ukasz Stettner$^{\ast}$}
\address{$^{\ast}$Institute of Mathematics, Polish Academy of Sciences, Warsaw, Poland}
\thanks{$^{\dagger}$Corresponding author}
\thanks{$^{\ast}$Damian Jelito and {\L}ukasz Stettner acknowledge research support by Polish National Science Centre grant no. 2020/37/B/ST1/00463. Part of this work was completed with the help of  University of Warsaw grant IDUB-
POB3-D110-003/2022 and Simons Semester ``Stochastic Modeling and Control'' at IMPAN in May-June 2023.}
\email{l.stettner@impan.pl}

\makeatletter
\def\namedlabel#1#2{\begingroup
    #2%
    \def\@currentlabel{#2}%
    \phantomsection\label{#1}\endgroup
}
\makeatother

\begin{document}

\begin{abstract}
In this paper, we investigate the effects of applying generalised (non-exponential) discounting on a long-run impulse control problem for a Feller-Markov process. We show that the optimal value of the discounted problem is the same as the optimal value of its undiscounted version. Next, we prove that an optimal strategy for the undiscounted discrete time functional is also optimal for the discrete-time discounted criterion and nearly optimal for the continuous-time discounted one. This shows that the discounted problem, being time-inconsistent in nature, admits a time-consistent solution. Also, instead of a complex time-dependent Bellman equation one may consider its simpler time-independent version.


\bigskip
\noindent \textbf{MSC2020 subject classifications:} 
93E20, 49J21, 49K21, 60J25 
\end{abstract}

\maketitle


\section{Introduction}
In recent years, control problems with generalised discounting have attracted considerable attention in the literature; see e.g.~\cite{HuaNgu2018, JasNow2020,BauJasNow2021} and references therein. This can be attributed to the observation that the classic exponential discount function improperly describes the behaviour of economic agents; we refer the reader to e.g.~\cite{LoewPre1992} for a detailed  discussion on this phenomenon and to~\cite{SalMos2011} for the statistical analysis rejecting the use of exponential discounting based on empirical data. On the other hand, control problems with non-exponential discounting usually pose serious technical problems. In the Markovian context, this is reflected in the fact that one needs to consider a non-linear time-inhomogeneous Bellman equation even if the underlying problem is time-homogeneous; see e.g.~\cite{BjoMur2014,BjoKhaMur2017,JasNow2020} for more details. In particular, typically the underlying problem is time-inconsistent and control equilibria are used instead of classical control strategies; see e.g.~\cite{HuaZho2021,BayZhaZho2021,BjoKhaMur2022}, and references therein.

In this paper, we combine the generalised discounting framework with the average cost per unit time impulse control setting. Impulse control provides a useful framework for introducing discrete-type interventions in continuous time dynamics which results in the extensive literature. In the financial context, impulse control setting could be used e.g. to model exchange rate control or portfolios with transaction costs; see e.g.~\cite{RunYas2018,BelChr2019}, and references therein for recent contributions. Other
applications of impulse control strategies include i.a. controlling epidemics,
ecosystems (optimal harvesting schemes), and finance (cash management), see e.g.~\cite{PiuPla2020,Erd2013,CorRod2012} and references. In the literature, long-run average cost impulse control was studied mostly in the undiscounted or classically discounted case; see e.g.~\cite{Rob1981,Ste1982} for some classical contributions and~\cite{PalSte2017,MenRob2018,Ste2022} for more recent ones. 

The main contribution of this paper is the proof that generalised discounting has no significant effect on the long-run impulse control problem. More specifically, we show that the optimal value of the discounted functional is equal to the optimal value of its undiscounted version; for details see Theorem~\ref{th:main_discrete} in the discrete-time setting and Theorem~\ref{th:main} in the continuous-time framework. Thus, despite the non-exponential form of the discount function, the discounted problem admits a time-consistent solution and instead of a complex time-dependent Bellman equation one may consider its simpler time-independent version. Moreover, an optimal strategy for the discrete-time undiscounted problem is also optimal for the discrete-time discounted functional and nearly optimal for the continuous-time discounted one; see Corollary~\ref{cor:nearly_optimal}. This facilitates the numerical approximation of an optimal strategy and shows that one may effectively restrict attention to the undiscounted discrete-time long-run avarage cost problem.

This paper continues the research initiated in~\cite{Ste2023}, where a generic (discrete time) Markov decision process with generalised discounting was considered. In that paper, it was shown that the optimal value of the long-run functional is unaffected when the non-linear discounting is introduced. In the present paper, we show a similar result in the long-run impulse control framework. Additionally, we obtain the analogous characterisation in the continuous time setting which was not studied in~\cite{Ste2023}. This suggests that similar effects should be observed in any long-run control problems provided that the (generalised) discount function satisfies suitable growth conditions; see the discussion following Assumption~\eqref{Discount} for details.

The structure of this paper is as follows. This introduction is followed by Section~\ref{S:preliminaries}, where we set up the notation, formally introduce the problems, and discuss the assumptions. Next, in Section~\ref{S:discrete} we focus on discrete-time (discrete-cost) problems. In Proposition~\ref{pr:existence_discounted} and Corollary~\ref{cor:existence_undiscounted} we establish the existence of solutions to the associated Bellman equations for the discounted and undiscounted case. This is used in Theorem~\ref{th:main_discrete}, where we show the equivalence between these two problems, i.e. the fact that the optimal values and the optimal strategies are the same. We extend this result in Theorem~\ref{th:same_payoff_disc} and show that the discounted and undiscounted functionals bring the same value for any (not necessarily optimal) stationary Markov strategies.  Next, in Section~\ref{S:continuous}, we generalise our results to the continuous-time framework; see Theorem~\ref{th:main} for the main contribution. The result follows from the series of approximation arguments that can be of independent interest; see Theorem~\ref{th:approx_disc}, Theorem~\ref{th:approx_undisc}, Theorem~\ref{th:approx_dyadic_disc}, and Theorem~\ref{th:approx_dyadic_undisc} for details. In particular, they justify some natural approximation schemes associated with discrete-time versions of the problems. Also, in this way we avoid technical difficulties associated with solving continuous time Bellman equations.

\section{Preliminaries}\label{S:preliminaries}

Let $X=(X_t)_{t\geq 0}$ be a continuous time standard Feller--Markov process on a filtered probability space $(\Omega, \mathcal{F}, (\mathcal{F}_t), \mathbb{P})$. The process $X$ takes values in a compact space $E$ endowed with a metric $\rho$ and the Borel $\sigma$-field $\mathcal{E}$. With any $x\in E$ we associate a probability measure $\mathbb{P}_x$ describing the dynamics of the process $X$ starting in $x$; see Section 1.4 in~\cite{Shi1978} for details. Also, we use $\mathbb{E}_x$, $x\in E$, and  $P_t(x,A):=\mathbb{P}_x[X_t\in A]$, $t\geq 0$, $x\in E$, $A\in \mathcal{E}$, for the corresponding expectation operator and the transition probability, respectively.  
Next, for any metric space $Z$, by $\mathcal{C}_b(Z)$ we denote the family of continuous bounded real-valued functions on $Z$. With this family we consider the supremum norm $\Vert f\Vert :=\sup_{x\in Z}|f(x)|$ and the span semi-norm $\Vert f\Vert_{sp}:=\sup_{x,y\in Z}(f(x)-f(y))$, $f\in \mathcal{C}_b(Z)$. Finally, we follow the convention that $\mathbb{N}$ consists of non-negative integers (including $0$), and we set $\sum_{i=0}^{-1}(\cdot):=0$ and $\inf\emptyset :=\infty$.

In this paper we assume that a decision-maker is allowed to shift the controlled process to any point in $E$, i.e. the set of admissible after-impulse states, which is denoted by $U$, satisfies  $U=E$. These shifts are accomplished with the help of an impulse control strategy, i.e. a sequence $V:=(\tau_i,\xi_i)_{i=1}^\infty$, where $(\tau_i)$ is an increasing sequence of stopping times and $(\xi_i)$ is a sequence of $\mathcal{F}_{\tau_i}$-measurable after-impulse states with values in $U$. With any starting point $x\in E$ and a strategy $V$ we associate a probability measure $\mathbb{P}_{x}^V$ and the corresponding expectation operator $\mathbb{E}_x^V$. Under this measure, the controlled process $Y$ process starts at $x$ and follows its usual (uncontrolled) dynamics up to the time $\tau_1$. Then, it is immediately shifted to $\xi_1$ and starts its evolution again, etc. More formally, we consider a countable product of filtered spaces $(\Omega, \mathcal{F}, (\mathcal{F}_t))$ and a coordinate process $(X_t^1, X_t^2, \ldots)$. Then, we define the controlled process $Y$ as $Y_t:=X_t^i$, $t\in [\tau_{i-1},\tau_i)$ with the convention $\tau_0\equiv 0$. Under the measure $\mathbb{P}_{x}^V$ we get $Y_{\tau_i}=\xi_i$; we refer to Chapter V in~\cite{Rob1978} for the construction details; see also Appendix in~\cite{Chr2014} and Section 2 in~\cite{Ste1982}. 
The family of all impulse control strategies is denoted by $\mathbb{V}$.  
Also, note that, to simplify the notation, by $Y_{\tau_i^-}:=X_{\tau_i}^i$, $i\in \mathbb{N}_{*}$, we denote the state of the process right before the $i$th impulse (yet, possibly, after the jump). Finally, note that in the discrete cost problem, we may impose an additional restriction on the after-impulse states and we assume that the set $U\subset E$ is compact; see Section~\ref{S:discrete} for details.

In this paper, we are interested in the properties of  long-run discounted and undiscounted impulse control problems. More specifically, for any  $x\in E$ and $V\in \mathbb{V}$, we consider the functionals
\begin{align}\label{eq:funct_time_inc}
    J^d(x,V)&:=\limsup_{T\to\infty}\frac{\mathbb{E}_{x}^V\left[\int_0^{T} \beta(s)g(Y_{s})ds+\sum_{i=1}^\infty 1_{\{\tau_i\leq T\}}\beta(\tau_i)c(X_{\tau_i}^i,X_{\tau_i}^{i+1}) \right]}{\int_0^{T} \beta(s)ds},\\
    J(x,V)&:=\limsup_{T\to\infty}\frac{\mathbb{E}_{x}^V\left[\int_0^T g(Y_s)ds+\sum_{i=1}^\infty 1_{\{\tau_i\leq T\}}c(X_{\tau_i}^i,X_{\tau_i}^{i+1}) \right]}{T},
    \label{eq:funct_time_cons}
\end{align}
where $\beta$, $g$, and $c$ are a discount function, a running reward/cost function, and a shift cost function, respectively. Note that~\eqref{eq:funct_time_cons} is the standard undiscounted average cost per unit time functional while~\eqref{eq:funct_time_inc} could be seen as its discounted version with discount function $\beta$. Also, $J$ corresponds to $J^d$ with $\beta\equiv 1$.

The aim of this paper is to characterise
\begin{equation}\label{eq:goal_cont}
    \inf_{V\in \mathbb{V}} J^d(x,V) \quad \text{and} \quad \inf_{V\in \mathbb{V}} J(x,V), \quad  x\in E.
\end{equation}
In fact, we show that the optimal values in~\eqref{eq:goal_cont} are the same; see Theorem~\ref{th:main} for details.

In the analysis, we assume the following conditions.
\begin{enumerate}
\item[(\namedlabel{Discount}{$\mathcal{A}1$})] (Discount function).  The map $\beta\colon [0,\infty)\ni t \mapsto \beta(t) \in (0,1]$ is (weakly) decreasing, continuous and satisfies $\beta(0)=1$. Also, we assume that $\beta$ is supermultiplicative, i.e. 
\begin{equation}\label{eq:beta_mult}
    \beta(t+s)\geq \beta(t)\beta(s), \quad t,s\geq 0,
\end{equation}
and we have $\int_0^\infty \beta(s)ds = \infty$. 
\item[(\namedlabel{A1}{$\mathcal{A}2$})] (Cost functions). The map $g:E\to \mathbb{R}$ is continuous and bounded. Also, the map $c:E\times U \to \mathbb{R}_{+}$ is  continuous, bounded,  strictly positive, and satisfies the triangle inequality, i.e. for some $c_0>0$, we have
\begin{equation}\label{eq:cost_ineq}
  0<c_0\leq c(x,\xi)\leq c(x,\eta)+c(\eta,\xi), \quad x,\xi,\eta\in E.  
\end{equation}
\item[(\namedlabel{Feller}{$\mathcal{A}3$})] (Feller property). The uncontrolled Markov process satisfies the $\mathcal{C}_b$-Feller property, i.e. for any $f\in \mathcal{C}_b(E)$ and $t>0$, the map $x\mapsto \mathbb{E}_x[f(X_t)]$ belongs to $\mathcal{C}_b(E)$.
\item[(\namedlabel{A2}{$\mathcal{A}4$})] (Uniform ergodicity). For any $t>0$ there exists $\Lambda_t\in (0,1)$ such that
\[
\sup_{x,y\in E}\sup_{A\in \mathcal{E}}|P_t(x,A)-P_t(y,A)|\leq \Lambda_t.
\]
\end{enumerate}

Let us now provide a more extensive comment on these conditions. 

Assumption~\eqref{Discount} states the properties of the considered discount function. In particular, monotonicity is a standard requirement to obtain a proper structure of the time value of money (the longer we wait for the cash inflow, the less valuable it is from today's perspective). Note that typically it is also assumed that $\lim_{t\to\infty}\beta(t)=0$. However, in our analysis, we do not directly use this assumption and we may cover both discounted and undiscounted ($\beta \equiv 1$) case. Next, note that the inequality stated in~\eqref{eq:beta_mult} could be linked to the decreasing impatience phenomenon observed in economic studies. In a nutshell, one may consider the behaviour of an individual facing two types of choices: (a) 100\$ today or 101\$ tomorrow, (b) 100\$ in 30 days or 101\$ in 31 days. It was observed that in case (a) people tend to choose an immediate reward of 100\$ while in (b) they are willing to wait a bit longer to get 101\$; see e.g.~\cite{LoewPre1992,Pre2004}, and references therein for details. This observation could be quantified by the ``forward-rate'' inequality $\beta(t+s)/\beta(t)\geq \beta(s)$, which is equivalent to~\eqref{eq:beta_mult}. Next, the property $\int_0^\infty \beta(s)ds=\infty$ states that the discount function does not vanish too fast. In fact, this is a crucial property used in our arguments as it facilitates suitable ergodic properties of the problems. More specifically, if $\int_0^\infty \beta(s)ds<\infty$, then the denominator in~\eqref{eq:funct_time_inc} converges to a finite constant and it is not expected that the associated optimisation problem is ergodic. In particular, the optimal value of the problem might depend on the initial state and a different form of the Bellman equation is needed; see Proposition~\ref{pr:opt_str_disc_disc} for the argument details and e.g.~\cite{Rob1981} for the analysis of the infinite horizon impulse control problem with (classical) discounting. Finally, it should be noted that all conditions from~\eqref{Discount} are satisfied e.g. for the generalised hyperbolic discount function $\beta(t):=\frac{1}{(1+ht)^{\alpha}}$, $t\geq 0$, with $h>0$ and $\alpha\in (0,1]$.

Assumption~\eqref{A1} states various properties of the running cost function and the shift-cost function. In particular,~\eqref{eq:cost_ineq} prevents multiple impulses at the same time. More specifically, when considering an impulse from $x$ to $\eta$ followed by an immediate impulse from $\eta$ to $\xi$, based on~\eqref{eq:cost_ineq}, a decision-maker should directly shift the process from $x$ to $\xi$. 
Typical examples satisfying~\eqref{A1} are in the form $c(x,\xi)=h(\rho(x,\xi))+c_0$, $x\in E$, $\xi\in U$, where $h$ is non-negative, bounded, increasing and subadditive, $c_0>0$, and $\rho$ is the underlying metric on $E$. Also, the assumption $0<c_0<c(x,\xi)$, $x,\xi\in E$, could be associated with some minimal (positive) fixed cost of any intervention. This assumption ensures i.a. that, for optimal strategies, there are finitely many impulses in the finite time interval; see the discussion following Equation~\ref{eq:Fatou} for details.

Assumption~\eqref{Feller} requires the $\mathcal{C}_b$-Feller property of the underlying Markov process and is a standard condition in the control literature. 

Assumption~\eqref{A2} requires uniform ergodicity of the underlying Markov process. This type of assumption is frequently used in the average cost framework to obtain the existence of a solution to the associated optimality equation; see e.g. Section 3.3 in~\cite{Her1989} and references therein for details. 

Finally, it should be noted that both Assumption~\eqref{Feller} and Assumption~\eqref{A2} are satisfied for regular reflected diffusions (with possible jumps) in a compact domain; see~\cite{MenRob1997} for details. In this case, the $\mathcal{C}_b$-Feller property follows from the existence of the continuous and bounded density (the Green function) for the transition semigroup. Also, Assumption~\eqref{A2} follows directly from Theorem 2.1 in~\cite{MenRob1997}. We refer to Section 7.1 in~\cite{GarMen2002} for a more detailed discussion.

\section{Discrete cost problem}\label{S:discrete}

This section considers a discretised version of the impulse control problems stated in~\eqref{eq:goal_cont}. More specifically, throughout this section we fix some $h>0$ and study the impulse control strategies with impulse times restricted to the time grid $\{0,h,2h, \ldots\}$ which affect discrete-type objective functionals. Also, in this section we may consider a more general control problem, i.e. we assume that $E$ is only a locally compact separable metric space and the after-impulse values can be restricted to some compact set $U\subset E$. Before we state explicit formulae for the control problems, let us introduce a  discrete-time discount function
\begin{equation}\label{eq:phi_h}
    \varphi_h(i):=\frac{1}{h}\int_{ih}^{(i+1)h}\beta(s)ds, \quad i\in \mathbb{N}.
\end{equation}
In the following lemma, we show that most of the properties of $\beta$ stated in~\eqref{Discount} transfers to $\varphi_h$.
\begin{lemma}\label{lm:discount_discrete}
Let $\varphi_h$ be given by~\eqref{eq:phi_h}. Then, we get $\varphi_h(0)\leq 1$ and the map $i\mapsto \varphi_h(i)$ is decreasing. Also, we get $\sum_{i=1}^\infty\varphi_h(i)=\infty$ and
\[
\varphi_h(i+k)\geq \varphi_h(i)\varphi_h(k), \quad i,k\in \mathbb{N}.
\]
\end{lemma}
\begin{proof}
Note that, using the monotonicity of $\beta$ (see Assumption~\eqref{Discount}), we have 
\[
\varphi_h(i+1)=\frac{1}{h}\int_{(i+1)h}^{(i+2)h}\beta(s)ds\leq\frac{1}{h}h\beta((i+1)h)\leq \frac{1}{h}\int_{ih}^{(i+1)h}\beta(s)ds=\varphi_h(i), \quad i \in \mathbb{N}.
\]
Thus, the map $i\mapsto \varphi_h(i)$ is decreasing. Using a similar argument we may show $\varphi(0)\leq 1$. Also, recalling that $\int_0^\infty \beta(s)ds=\infty$, we get $\sum_{i=1}^\infty\varphi_h(i) = \frac{1}{h}\int_0^\infty \beta(s)ds=\infty$. Next, let us fix some $i,k\in \mathbb{N}$. Then, we get
\begin{align*}
    \varphi_h(i)\varphi_h(k) &= \frac{1}{h^2}\int_{ih}^{(i+1)h}\int_{kh}^{(k+1)h}\beta(s) \beta(u) du ds \\
    &\leq \frac{1}{h^2}\int_{ih}^{(i+1)h}\int_{kh}^{(k+1)h}\beta(s+u)  du ds\\
    & = \frac{1}{h^2}\int_{ih}^{(i+1)h}\int_{v_1+kh}^{v_1+(k+1)h}\beta(v_2)  dv_2 dv_1\\
    & \leq \frac{1}{h^2}\int_{ih}^{(i+1)h}\beta(v_1+kh)h dv_1    = \frac{1}{h}\int_{(i+k)h}^{(i+k+1)h}\beta(s)  ds = \varphi_h(i+k),
\end{align*}
which concludes the proof.
\end{proof}

Next, by $\mathbb{V}_h\subset \mathbb{V}$  we denote the subfamily of impulse control strategies with impulse times restricted to $\{0,h,2h, \ldots\}$. Also, for any  $x\in E$ and $V\in \mathbb{V}_h$, we define the discretised versions of the functionals~\eqref{eq:funct_time_inc} and~\eqref{eq:funct_time_cons} by the formulae
\begin{align}\label{eq:funct_time_inc_disc}
    J^d_h(x,V)&:=\limsup_{n\to\infty}\frac{1}{\sum_{i=0}^{n-1}\varphi_h(i)}\mathbb{E}_{x}^V\left[\sum_{i=0}^{n-1} \varphi_h(i)g(Y_{ih})\right.\nonumber\\
    &\phantom{=}\left.+\sum_{i=1}^\infty 1_{\{\tau_i\leq nh\}}\frac{1}{h}\varphi_h(\tau_i/h)c(X_{\tau_i}^i,X_{\tau_i}^{i+1}) \right],\\
    \label{eq:funct_time_cons_disc}
    J_h(x,V)&:=\limsup_{n\to\infty}\frac{1}{n}\mathbb{E}_{x}^V\left[\sum_{i=1}^n g(Y_{ih})+\sum_{i=1}^\infty 1_{\{\tau_i\leq nh\}}\frac{1}{h}c(X_{\tau_i}^i,X_{\tau_i}^{i+1}) \right];
\end{align}
note that here the normalisation by $h$ ensures that these functionals converge to $J^d$ and $J$ as $h\to 0$; see Theorem~\ref{th:approx_disc} and Theorem~\ref{th:approx_dyadic_undisc} for details.

By analogy to~\eqref{eq:goal_cont}, the main goal of this section is to characterise
\begin{align}\label{eq:goal_discrete}
    \inf_{V\in \mathbb{V}_{h}} J^d_h(x,V), \quad \text{and} \quad     \inf_{V\in \mathbb{V}_h}  J_h(x,V), \quad x\in E.
\end{align}
To do this, we prove the existence of solutions $w^d_h\in \mathcal{C}_b(\mathbb{N}\times E)$, $w_h\in \mathcal{C}_b(E)$, $\lambda^d_h(k)\in \mathbb{R}$ for $k\in \mathbb{N}$, and $\lambda_h\in \mathbb{R}$ to two types of Bellman equations
\begin{align}
w^d_h(k,x)&=\min\left(\varphi_h(k)\left(g(x)-\lambda^d_h(k)\right)+\mathbb{E}_{x}\left[w^d_h(k+1,X_{h})\right],M^d_h w^d_h(k,x)\right),\label{eq:Bellman_disc_inc}\\
    w_h(x)& = \min\left(g(x)-\lambda_h+\mathbb{E}_x[w_h(X_h)], M_h w_h(x)\right),\label{eq:Bellman_disc_cons}
\end{align}
where $k\in \mathbb{N}$, $x\in E$, and the operators $M^d_h$ and $M_h$ are given as
\begin{align}\label{eq:M^d_h}
    M^d_h w^d_h(k,x)&:=\inf_{\xi\in U}\left(\varphi_h(k)\frac{1}{h}c(x,\xi)+w^d_h(k,\xi)\right), \quad k\in \mathbb{N}, \, x\in E,\\
M_h w_h(x)&:=\inf_{\xi\in U}\left(\frac{1}{h}c(x,\xi)+w_h(\xi)\right), \quad  x\in E.\label{eq:M_h}
\end{align}
Note that the links between~\eqref{eq:Bellman_disc_inc} and~\eqref{eq:Bellman_disc_cons} and the optimal values defined in~\eqref{eq:goal_discrete} are shown in Proposition~\ref{pr:opt_str_disc_disc} and Corollary~\ref{cor:opt_str_disc_undisc}. Also, note that~\eqref{eq:Bellman_disc_inc} could be seen as a time-dependent version of~\eqref{eq:Bellman_disc_cons}, and this time-dependence follows from the fact that we consider a generic, non-constant discount function.

We start with establishing the existence of solutions to~\eqref{eq:Bellman_disc_inc} and~\eqref{eq:Bellman_disc_cons}. In the proof, we use the properties of the equivalence relation $\sim$ defined for $f,g\in \mathcal{C}_b(\mathbb{N}\times E)$ as $f\sim g$ if and only if $\sup_{i\in \mathbb{N}}\Vert f(i,\cdot)-g(i,\cdot)\Vert=0$, where $\Vert \cdot \Vert$ is the supremum norm (on $E$). With this relation we define the quotient space $\mathcal{C}_{sp}:=\mathcal{C}_b(\mathbb{N}\times E)/\sim$ and we equip it with the the norm $\vvvert f\vvvert:=\sup_{i\in \mathbb{N}}\Vert f(i,\cdot)\Vert_{sp}$. It can be shown that $\mathcal{C}_{sp}$ with $\vvvert \cdot \vvvert$ is a Banach space; see e.g. Definition 3.4 and Lemma 3.5 in Chapter 3 of~\cite{Her1989} for further discussion.

\begin{proposition}\label{pr:existence_discounted}
    There exists a unique up to a constant map $w^d_h\in \mathcal{C}_b(\mathbb{N}\times E)$ and a unique sequence of constants $(\lambda^d_h(k))_{k\in \mathbb{N}}\subset \mathbb{R}$ which are a solution to Equation~\eqref{eq:Bellman_disc_inc}.
\end{proposition}
\begin{proof}
In the proof we set $h=1$; the general case follows the same logic. Also, we show that there exist $w^d_1\in \mathcal{C}_b(\mathbb{N}\times E)$ and $(\lambda^d_h(k))_{k\in \mathbb{N}}\subset \mathbb{R}$ satisfying, for any $k\in \mathbb{N}$ and $x\in E$, the relation
\begin{align}\label{eq:Bellman_disc_inc_existence}
    w^d_1(k,x)&=\min\bigg(\varphi_1(k)\left(g(x)-\lambda^d_1(k)\right)+\mathbb{E}_{x}\left[w^d_1(k+1,X_{1})\right],\nonumber \\
&\phantom{=====}\inf_{\xi\in U} \left(\varphi_1(k)c(x,\xi)+\varphi_1(k)\left(g(\xi)-\lambda^d_1(k)\right)+\mathbb{E}_{\xi}\left[w^d_1(k+1,X_{1})\right]\right)\bigg).
\end{align}
Then, using the argument from Theorem 3.1 in~\cite{JelPitSte2019b}, we get that the pair $(w^d_1,(\lambda^d_1(k)))$ solves~\eqref{eq:Bellman_disc_inc}.

To solve~\eqref{eq:Bellman_disc_inc_existence}, we use the span contraction approach. To do this, we set $\bar{c}(x,\xi):=c(x,\xi) 1_{\{x\neq \xi\}}$, $x\in E$, $\xi\in U$. Also, we define the operator 
\begin{align*}
    Fv(k,x)&:=\inf_{\xi\in U\cup \{x\}}\left(\varphi_1(k)\bar{c}(x,\xi)+\varphi_1(k)g(\xi)+\mathbb{E}_{\xi}\left[v(k+1,X_{1})\right]\right),
\end{align*}
where $v\in \mathcal{C}_b(\mathbb{N}\times E)$, $k\in \mathbb{N}$, $x\in E$; note that the operator, up to the constant, is simply given by the right-hand side of~\eqref{eq:Bellman_disc_inc_existence}. In particular $F\mathcal{C}_b(\mathbb{N}\times E)\subset \mathcal{C}_b(\mathbb{N}\times E)$ and Equation~\eqref{eq:Bellman_disc_inc_existence} is equivalent to
\begin{equation}\label{eq:Bellman_disc_inc_fixed_point}
    w^d_1(k,x)+\varphi_1(k)\lambda^d_1(k) = Fw^d_1(k,x), \quad k\in \mathbb{N}, \, x\in E.
\end{equation}

We now show that $F$ is a contraction in the norm $\vvvert\cdot \vvvert$ on $\mathcal{C}_{sp}$. Let $v_1,v_2\in \mathcal{C}_b(\mathbb{N}\times E)$, $k\in \mathbb{N}$, $x,y\in E$. Also, let $\xi_1$ and $\xi_2$ be minimisers for $Fv_2(k,x)$ and $Fv_1(k,y)$, respectively. Then, using~\eqref{A2}, we get
\begin{align}\label{eq:pr:existence_discounted:span}
    Fv_1&(k,x)-Fv_2(k,x)-Fv_1(k,y)+Fv_2(k,y)\nonumber\\
    & \leq \varphi_1(k)\bar{c}(x,\xi_1)+\varphi_1(k)g(\xi_1)+\mathbb{E}_{\xi_1}\left[v_1(k+1,X_{1})\right]\nonumber\\
    & \phantom{=} -  \varphi_1(k)\bar{c}(x,\xi_1)-\varphi_1(k)g(\xi_1)-\mathbb{E}_{\xi_1}\left[v_2(k+1,X_{1})\right]\nonumber\\
    & \phantom{=} -  \varphi_1(k)\bar{c}(x,\xi_2)-\varphi_1(k)g(\xi_2)-\mathbb{E}_{\xi_2}\left[v_1(k+1,X_{1})\right]\nonumber\\
    & \phantom{=} + \varphi_1(k)\bar{c}(x,\xi_2)+\varphi_1(k)g(\xi_2)+\mathbb{E}_{\xi_2}\left[v_2(k+1,X_{1})\right]\nonumber\\
    &=\int_E(v_1(k+1,z_1)-v_2(k+1,z_1))P_1(\xi_1,dz_1)\nonumber\\
    &\phantom{=}-\int_E(v_1(k+1,z_2) -v_2(k+1,z_2))P_1(\xi_2,dz_2)\nonumber\\
    &\leq \Lambda_1\Vert v_1(k+1,\cdot)-v_2(k+1,\cdot)\Vert_{sp} .
\end{align}
Consequently, we get
\[
\vvvert Fv_1-Fv_2\vvvert\leq \Lambda_1\vvvert v_1-v_2\vvvert
\]
and the operator $F$ is a contraction with respect to $\vvvert \cdot \vvvert$. Thus, using the Banach fixed-point theorem, we get that there exists a map $w^d_1\in \mathcal{C}_{sp}$ such that $\vvvert Fw^d_1-w^d_1 \vvvert=0$. This implies that, for any $k\in \mathbb{N}$, there exists $\lambda_1^d(k)\in \mathbb{R}$ such that~\eqref{eq:Bellman_disc_inc_fixed_point} is satisfied. Moreover, directly from the definition of the span norm, we get the map $w^d_1$ is unique up to a constant and $\lambda_1^d(k)$ is unique, which concludes the proof.
\end{proof}

Using the argument from Proposition~\ref{pr:existence_discounted} we may also show the existence of a solution to the undiscounted Bellman equation~\eqref{eq:Bellman_disc_cons}; see Corollary~\ref{cor:existence_undiscounted}. In fact, it is enough to set $\varphi_h\equiv 1$ which corresponds to $\beta \equiv 1$. 

\begin{corollary}\label{cor:existence_undiscounted}
There exists a unique up to a constant map $w_h\in \mathcal{C}_b(E)$ and a unique constant $\lambda_h\in \mathbb{R}$ which are a solution to Equation~\eqref{eq:Bellman_disc_cons}.
\end{corollary}

Next, we provide a martingale characterisation of~\eqref{eq:Bellman_disc_inc}.

\begin{proposition}\label{pr:martingale_disc_inc}
Let $(w^d_h,\lambda^d_h)$ be a solution to~\eqref{eq:Bellman_disc_inc}. Also, for any $k\in \mathbb{N}$, let us define
\begin{align*}
    z_k^d(n)&:=\sum_{i=k}^{n+k-1} \varphi_h(i)(g(X_{(i-k)h})-\lambda^d_h(i))+w(k+n,X_{nh}), \quad n\in \mathbb{N},\\
    \hat\tau_k^d&:=h\inf\{n\in \mathbb{N}\colon w^d_h(k+n,X_{nh})=M^d_hw^d_h(k+n,X_{nh})\}.
\end{align*}
Then, for any $k\in \mathbb{N}$ and $x\in E$, the process $(z_k^d(n))_{n\in \mathbb{N}}$ is a $\mathbb{P}_x$-submartingale and the process $(z_k^d(n\wedge (\hat\tau_k^d/h)))_{n\in \mathbb{N}}$ is a $\mathbb{P}_x$-martingale.
\end{proposition}
\begin{proof}
For transparency, we show the proof only for $h=1$; the general case follows the same logic. Let $x\in E$ and $k,n\in \mathbb{N}$. First, we show the submartingale property. Recalling~\eqref{eq:Bellman_disc_inc}, we get
\[
\varphi_1(n+k)(g(X_{n})-\lambda^d(n+k))+\mathbb{E}_{X_{n}}[w_1(n+k+1,X_1)]\geq w(n+k,X_{n}) \quad \mathbb{P}_x \, \text{a.s.}
\]
Thus, we have
\begin{align*}
    \mathbb{E}_x[z_k^d(n+1)|\mathcal{F}_n] & = \sum_{i=k}^{n+k-1} \varphi_h(i)(g(X_{i-k})-\lambda^d(i))+\varphi_1(n+k)(g(X_{n})-\lambda^d(n+k))\\
    & \phantom{=}+\mathbb{E}_{X_{n}}[w_1(n+k+1,X_1)] \\
    & \geq \sum_{i=k}^{n+k-1} \varphi_h(i)(g(X_{i-k})-\lambda^d(i))+w(n+k,X_{n}) = z_k^d(n),
\end{align*}
which shows that $(z_k^d(n))_{n\in \mathbb{N}}$ is a submartingale.

Second, we show the martingale property. Using~\eqref{eq:Bellman_disc_inc} and the definition of $\hat\tau^d_1$, on the event $\{\hat\tau^d_1\geq n+1\}$, we get
 \[
\varphi_1(n+k)(g(X_{n})-\lambda^d(n+k))+\mathbb{E}_{X_{n}}[w_1(n+k+1,X_1)]= w(n+k,X_{n}) \quad \mathbb{P}_x \, \text{a.s.}
\]
Thus, we obtain
\begin{align*}
    \mathbb{E}_x&[z_k^d((n+1)\wedge\hat\tau^d_1)|\mathcal{F}_n]  = \sum_{i=k}^{n\wedge \hat\tau^d_1+k-1} \varphi_h(i)(g(X_{i-k})-\lambda^d(i))+1_{\{\hat\tau^d_1\leq n\}} w(n+k,X_{n}) \\
    \\
    & \phantom{=}+1_{\{\hat\tau^d_1\geq n+1\}} \left(\varphi_1(n+k)(g(X_{n})-\lambda^d(n+k)) + \mathbb{E}_{X_{n}}[w_1(n+k+1,X_1)]\right) \\
    & =\sum_{i=k}^{n\wedge \hat\tau^d_1+k-1} \varphi_h(i)(g(X_{i-k})-\lambda^d(i))+w(n\wedge \hat\tau^d_1+k,X_{n\wedge \hat\tau^d_1}) = z_k^d(n\wedge \hat\tau^d_1),
\end{align*}
which concludes the proof.
 \end{proof}

Using the argument from Proposition~\ref{pr:martingale_disc_inc} with $\beta\equiv 1$ we also get a martingale characterisation of~\eqref{eq:Bellman_disc_cons}; see Corollary~\ref{cor:martingale_disc_cons}. 

\begin{corollary}\label{cor:martingale_disc_cons}
Let $(w_h,\lambda_h)$ be a solution to~\eqref{eq:Bellman_disc_cons}. Also, for any $k\in \mathbb{N}$, let us define
\begin{align*}
    z(n)&:=\sum_{i=0}^{n-1}(g(X_{ih})-\lambda_h)+w_h(X_{nh}), \quad n\in \mathbb{N},\\
    \hat\tau&:=h\inf\{n\in \mathbb{N}\colon w_h(X_{nh})=M_h w_h(X_{nh})\}.
\end{align*}
Then, for any $x\in E$, the process $(z(n))_{n\in \mathbb{N}}$ is a $\mathbb{P}_x$-submartingale and the process $(z(n\wedge (\hat\tau/h)))_{n\in \mathbb{N}}$ is a $\mathbb{P}_x$-martingale.
\end{corollary}

Now, we show how to link~\eqref{eq:Bellman_disc_inc} and~\eqref{eq:Bellman_disc_cons} with optimal strategies for~\eqref{eq:goal_discrete}. We start with introducing some notation.  Let $(w^d_h,\lambda^d_h)$ be a solution to~\eqref{eq:Bellman_disc_inc}. With this pair and any $k\in \mathbb{N}$, we associate a strategy $\hat{V}^d_h:=(\hat\sigma^d_{i},\hat\xi^d_i)_{i=1}^\infty$ given recursively by
\begin{align}\label{eq:opt_str_disc_inc}
    \hat\sigma^d_{i}&:=h\inf\{n\in \mathbb{N}\colon n\geq \hat\sigma^d_{i-1}/h, \, w^d_h(n,X_{nh}^{i})=M^d_hw^d_h(n,X_{nh}^i)\}, \nonumber\\
    \hat\xi^d_i&:=\argmin_{\xi\in U}(\varphi_h(\hat\sigma^d_{i})\frac{1}{h}c(X^i_{\hat\sigma^d_{i}},\xi)+w^d_h(\hat\sigma^d_{i},\xi))1_{\{\hat\sigma^d_{i}<\infty\}}+\xi_0 1_{\{\hat\sigma^d_{i}=\infty\}},
\end{align}
where we set $\hat\sigma_0^d\equiv 0$, and $\xi_0\in U$ is some fixed point. Similarly, for $(w_h,\lambda_h)$ being a solution to~\eqref{eq:Bellman_disc_cons}, we define a strategy $\hat{V}_h:=(\hat\sigma_i,\hat\xi_i)_{i=1}^\infty$ by
\begin{align}\label{eq:opt_str_disc_cons}
    \hat\sigma_i&:=h\inf\{n\in \mathbb{N}\colon n\geq \hat\sigma_{i-1}/h, \, w_h(X_{nh}^{i})=M_h w_h(X_{nh}^i)\}, \nonumber\\
    \hat\xi_i&:=\argmin_{\xi\in U}(\frac{1}{h}c(X^i_{\hat\sigma_{i}},\xi)+w_h(\xi)) 1_{\{\hat\sigma_{i}<\infty\}}+\xi_0 1_{\{\hat\sigma_{i}=\infty\}},
\end{align}
where $\hat\sigma_0\equiv 0$ and $\xi_0\in U$ is some fixed point. Note that the after-impulse states $\hat\xi_i^d$ and $\hat\xi_i$ are given by the minimisers of the operators $M^d_h$ and $M_h$ from~\eqref{eq:M^d_h} and~\eqref{eq:M_h}. Also, the impulse times $\hat\sigma_i^d$ and $\hat\sigma_i$ are linked to the stopping times from Proposition~\ref{pr:martingale_disc_inc} and Corollary~\ref{cor:martingale_disc_cons}. Finally, recall from Section~\ref{S:preliminaries} that we follow the convention $\inf\emptyset:=\infty$.

\begin{proposition}\label{pr:opt_str_disc_disc}
Let $(w^d_h,\lambda^d_h)$ be a solution to~\eqref{eq:Bellman_disc_inc}. Also, let the strategy $\hat{V}^d_h$ be given by~\eqref{eq:opt_str_disc_inc}, and let $J^d_h$ be given by~\eqref{eq:funct_time_inc_disc}. Then, we get
\[
\inf_{V\in \mathbb{V}_{h}} J^d_h(x,V)= J^d_h(x,\hat{V}^d_h)=\limsup_{n\to\infty}\frac{\sum_{i=0}^{n-1}\varphi_h(i)\lambda^d_h(i)}{\sum_{i=0}^{n-1}\varphi_h(i)}, \quad x\in E.
\]
\end{proposition}
\begin{proof}
    For the notational convenience, we set $h=1$; the general case follows the same logic. 

    First, we show that 
  \begin{equation}\label{eq:pr:ver_th_disc_inc:0}
      J^d_1(x,\hat{V}^d_1)=\limsup_{n\to\infty}\frac{\sum_{i=0}^{n-1}\varphi_1(i)\lambda^d_1(i)}{\sum_{i=0}^{n-1}\varphi_1(i)}, \quad x\in E.    
  \end{equation}
Using Proposition~\ref{pr:martingale_disc_inc}, we get that the process $z_1^d(n\wedge \hat\sigma_1^d)_{n\in \mathbb{N}}$ is a martingale. Also, note that on the event $\{\hat\sigma_1^d<\infty\}$, we have
\[
w^d_1(\hat\sigma_1^d, X_{\hat\sigma_1^d}^1)=M^d_1w^d_1(\hat\sigma_1^d, X_{\hat\sigma_1^d}^1)=\varphi_1(\hat\sigma_1^d)c(X_{\hat\sigma_1^d}^1,X_{\hat\sigma_1^d}^2)+w^d_1(\hat\sigma_1^d,X_{\hat\sigma_1^d}^2).
\]
Thus, for any $x\in E$ and $n\in \mathbb{N}$, we obtain
\begin{align*}
    w^d_1(0,x)&=\mathbb{E}_x^{\hat{V}^d_1}\left[\sum_{i=0}^{n\wedge \hat\sigma_1^d-1}\varphi_1(i)(g(X_i^1)-\lambda_1^d(i))+w_1^d(n\wedge\hat\sigma_1^d, X^1_{n\wedge\hat\sigma_1^d}) \right]\\
    & = \mathbb{E}_x^{\hat{V}^d_1}\left[\sum_{i=0}^{n\wedge \hat\sigma_1^d-1}\varphi_1(i)(g(X_i^1)-\lambda_1^d(i)) +1_{\{\hat\sigma_1^d\leq n \}}\varphi_1(\hat\sigma_1^d)c(X_{\hat\sigma_1^d}^1,X_{\hat\sigma_1^d}^2) \right.\\
    &\phantom{=\mathbb{E}_x^{\hat{V}^d}}\left.+1_{\{\hat\sigma_1^d\leq n \}}w^d_1(\hat\sigma_1^d, X^2_{\hat\sigma_1^d}) + 1_{\{\hat\sigma_1^d> n \}}w^d_1(n, X^1_{n})\right].
\end{align*}
Consequently, for any $m\in  \mathbb{N}$, recursively we get
\begin{align}\label{eq:pr:ver_th_disc_inc:1}
    w^d_1(0,x) & = \mathbb{E}_x^{\hat{V}^d_1}\left[\sum_{i=0}^{n\wedge \hat\sigma_m^d-1}\varphi_1(i)(g(Y_i)-\lambda_1^d(i)) +\sum_{i=1}^m1_{\{\hat\sigma_i^d\leq n \}}\varphi_1(\hat\sigma_i^d)c(X_{\hat\sigma_i^d}^i,X_{\hat\sigma_i^d}^{i+1}) \right.\nonumber\\
    &\phantom{=\mathbb{E}_x^{\hat{V}^d_1}}\left.+1_{\{\hat\sigma_m^d\leq n \}}w^d_1(\hat\sigma_m^d, X^m_{\hat\sigma_m^d}) + \sum_{i=1}^m 1_{\{\hat\sigma_{i-1}^d\leq n< \hat\sigma_i^d \}}w^d_1(n, X^i_{n})\right].
\end{align}
Hence, recalling that $0<c_0\leq c(\cdot,\cdot)$ (see Assumption~\eqref{A1}) and using the boundedness of $w^d_1$, we obtain
\begin{align*}
     \mathbb{E}_x^{\hat{V}^d_1}\left[\sum_{i=1}^m1_{\{\hat\sigma_m^d\leq n \}}\varphi_1(\hat\sigma_i^d)c_0 \right]&\leq \mathbb{E}_x^{\hat{V}^d_1}\left[\sum_{i=1}^m 1_{\{\hat\sigma_i^d\leq n \}}\varphi_1(\hat\sigma_i^d)c(X_{\hat\sigma_i^d}^i,X_{\hat\sigma_i^d}^{i+1}) \right]\\
     & \leq 2 \Vert w_1^d\Vert +n\Vert g\Vert +\sup_{i=0,\ldots, n}\lambda_1^d(i)<\infty
\end{align*}
Consequently, using Fatou's lemma, for any $n\in \mathbb{N}$, we get 
\begin{equation}\label{eq:Fatou}
    \mathbb{E}_x^{\hat{V}^d_1}\left[\sum_{i=1}^\infty 1_{\{\hat\sigma_m^d\leq n \}}\varphi_1(\hat\sigma_i^d)c_0 \right]<\infty,
\end{equation}
which, combined with the fact $c_0>0$, implies $\lim_{m\to\infty}\hat\sigma_m^d=\infty $. Thus, letting $m\to\infty$ in~\eqref{eq:pr:ver_th_disc_inc:1} and using Lebesgue’s dominated convergence theorem, we obtain
\begin{align}\label{eq:pr:ver_th_disc_inc:2}
    w^d_1(0,x) & = \mathbb{E}_x^{\hat{V}^d_1}\left[\sum_{i=0}^{n-1}\varphi_1(i)(g(Y_i)-\lambda_1^d(i)) +\sum_{i=1}^\infty 1_{\{\hat\sigma_i^d\leq n \}}\varphi_1(\hat\sigma_i^d)c(X_{\hat\sigma_i^d}^i,X_{\hat\sigma_i^d}^{i+1}) \right.\nonumber\\
    &\phantom{=\mathbb{E}_x^{\hat{V}^d}}\left. + \sum_{i=1}^\infty 1_{\{\hat\sigma_{i-1}^d\leq n< \hat\sigma_i^d \}}w^d_1(n, X^i_{n})\right].
\end{align}
Next, dividing both hand sides by $\sum_{i=0}^{n-1}\varphi_1(i)$ and using the boundedness of $w^d_1$, we have
\begin{multline}\label{eq:pr:ver_th_disc_inc:2.5}
    \frac{-2\Vert w_1^d \Vert}{\sum_{i=0}^{n-1}\varphi_1(i)}+\frac{\sum_{i=0}^{n-1}\varphi_1(i)\lambda_1^d(i))}{\sum_{i=0}^{n-1}\varphi_1(i)}\\
    \leq \frac{\mathbb{E}_x^{\hat{V}^d_1}\left[\sum_{i=0}^{n-1}\varphi_1(i)g(Y_i)+\sum_{i=1}^\infty 1_{\{\hat\sigma_i^d\leq n \}}\varphi_1(\hat\sigma_i^d)c(X_{\hat\sigma_i^d}^i,X_{\hat\sigma_i^d}^{i+1}) \right]}{\sum_{i=0}^{n-1}\varphi_1(i)}\\
    \leq \frac{2\Vert w_1^d \Vert}{\sum_{i=0}^{n-1}\varphi_1(i)}+\frac{\sum_{i=0}^{n-1}\varphi_1(i)\lambda_1^d(i))}{\sum_{i=0}^{n-1}\varphi_1(i)}.
\end{multline}
Also, letting $n\to\infty$ and recalling that by Lemma~\ref{lm:discount_discrete}, we have $\sum_{i=0}^{\infty}\varphi_1(i)=\infty$, we conclude that~\eqref{eq:pr:ver_th_disc_inc:0} is satisfied.

Second, we show that for any strategy $V:=(\tau_i,\xi)_{i=1}^\infty\in \mathbb{V}_1$, we have 
\begin{equation}\label{eq:pr:ver_th_disc_inc:3}
  \limsup_{n\to\infty}\frac{\sum_{i=0}^{n-1}\varphi_h(i)\lambda^d_h(i)}{\sum_{i=0}^{n-1}\varphi_h(i)}\leq J^d_h(x,V), \quad x\in E.   
\end{equation}
We start by noting that, from Proposition~\ref{pr:martingale_disc_inc} and Doob’s optional stopping theorem, the process $(z_1^d(n\wedge \tau_1))_{n\in \mathbb{N}}$ is a submartingale. Also, note that directly from~\eqref{eq:Bellman_disc_inc}, on the event $\{\tau_1<\infty\}$ we get
\[
w_1^d(\tau_1,X^1_{\tau_1})\leq M_1^dw_1^d(\tau_1,X^1_{\tau_1})\leq \varphi_1(\tau_1)c(X_{\tau_1}^1,\xi_1)+w_1^d(\tau_1,\xi_1).
\]
Thus, as in~\eqref{eq:pr:ver_th_disc_inc:2}, we obtain
\begin{align*}
    w^d_1(0,x) & \leq \mathbb{E}_x^{V}\left[\sum_{i=0}^{n-1}\varphi_1(i)(g(Y_i)-\lambda_1^d(i)) +\sum_{i=1}^\infty 1_{\{\tau_i\leq n \}}\varphi_1(\tau_i)c(X_{\tau_i}^i,X_{\tau_i}^{i+1}) \right.\nonumber\\
    &\phantom{=\mathbb{E}_x^{V}}\left. + \sum_{i=1}^\infty 1_{\{\tau_{i-1}^d\leq n< \tau_i \}}w^d_1(n, X^i_{n})\right].
\end{align*}
Hence, dividing by $\sum_{i=0}^{n-1}\varphi_1(i)$ and letting $n\to\infty$, we get~\eqref{eq:pr:ver_th_disc_inc:3}. Combining this with~\eqref{eq:pr:ver_th_disc_inc:0}, we conclude the proof.
\end{proof}

As in Proposition~\ref{pr:opt_str_disc_disc}, we may show that the strategy defined in~\eqref{eq:opt_str_disc_cons} is optimal for the functional~\eqref{eq:funct_time_cons_disc}. The proof is omitted for brevity.

\begin{corollary}\label{cor:opt_str_disc_undisc}
Let $(w_h,\lambda_h)$ be a solution to~\eqref{eq:Bellman_disc_cons}. Also, let the strategy $\hat{V}_h$ be given by~\eqref{eq:opt_str_disc_cons}, and let $J_h$ be given by~\eqref{eq:funct_time_cons_disc}. Then, we get
\[
\inf_{V\in \mathbb{V}_h} J_h(x,V)=J_h(x,\hat{V}_h)=\lambda_h, \quad x\in E.
\]
\end{corollary}

Now, we show the main result of this section, i.e. Theorem~\ref{th:main_discrete}. More specifically, we show that the strategy given by~\eqref{eq:opt_str_disc_cons}, which is optimal for the undiscounted functional~\eqref{eq:funct_time_cons_disc}, is also optimal for the discounted functional~\eqref{eq:funct_time_inc_disc}. Also, the optimal values defined in~\eqref{eq:goal_discrete} are the same. The argument is partially based on Theorem 2 from~\cite{Ste2023}.

\begin{theorem}\label{th:main_discrete}
Let $(w_h,\lambda_h)$ be a solution to~\eqref{eq:Bellman_disc_cons} and let the strategy $\hat{V}_h$ be given by~\eqref{eq:opt_str_disc_cons}. Also, let $J^d_h$ and $J_h$ be given by~\eqref{eq:funct_time_inc_disc} and~\eqref{eq:funct_time_cons_disc}, respectively. Then, we get
\[
 \inf_{V\in \mathbb{V}_h} J^d_h(x,V) =J^d_h(x,\hat{V}_h) =\lambda_h =\inf_{V\in \mathbb{V}_h}J_h(x,v), \quad x\in E.
\]
\end{theorem}
\begin{proof}
For notational convenience, we show the proof only for $h=1$; the general case follows the same logic. 

First, we show that for the strategy $\hat{V}_1$ given by~\eqref{eq:opt_str_disc_cons}, we get $J^d_1(x,\hat{V}_1)= \lambda_1$. Recalling~\eqref{eq:Bellman_disc_cons}, for any $k\in \mathbb{N}$ and $x\in E$, we get
\[
\varphi_1 (k) w_1(x)= \min\left(\varphi_1 (k) (g(x)-\lambda_1)+\varphi_1 (k) \mathbb{E}_x[w_1(X_1)], \varphi_1 (k) Mw_1(x)\right).
\]
Thus, for any $x\in E$, $n\in \mathbb{N}$, and $i=0, 1, \ldots, \hat\sigma_1 \wedge n-1$, we have
\begin{equation}\label{eq:th_main_disc:1}
  \varphi_1 (i) w_1(X_i^1)= \varphi_1 (i) (g(X_i^1)-\lambda_1)+\varphi_1 (i) \mathbb{E}_{X_i^1}^{\hat{V}_1}[w_1(X_1)], \quad \mathbb{P}_{x}^{\hat{V}_1} \, \text{a.s.}  
\end{equation}
Hence, using the fact that the controlled process is Markov between the consecutive impulses, we obtain
\begin{align*}
     0 & =\mathbb{E}_{x}^{\hat{V}_1}\left[\sum_{i=0}^{\hat\sigma_1 \wedge n-1}\varphi_1 (i) (g(X_i^1)-\lambda_1)+ \sum_{i=0}^{\hat\sigma_1 \wedge n-1}\varphi_1 (i) (\mathbb{E}_{X_i^1}^{\hat{V}_1}[w_1(X_1^1)]-w_1(X_i^1))\right]\\
     &=\mathbb{E}_{x}^{\hat{V}_1}\left[\sum_{i=0}^{\hat\sigma_1 \wedge n-1}\varphi_1 (i) (g(X_i^1)-\lambda_1)+ \sum_{i=0}^{\hat\sigma_1 \wedge n-1}\varphi_1 (i) (w_1(X_{i+1}^1)-w_1(X_i^1))\right].
\end{align*}

Next, note that on the event $\{\hat\sigma_1<\infty\}$, we have
\begin{align}\label{eq:th_main_disc:2}
\varphi_1(\hat\sigma_1)w_1(X_{\hat\sigma_1}^1) &= \varphi_1(\hat\sigma_1)Mw_1(X_{\hat\sigma_1}^1) = \varphi_1(\hat\sigma_1)c(X_{\hat\sigma_1}^1,X_{\hat\sigma_1}^2)+\varphi_1(\hat\sigma_1)w_1(X_{\hat\sigma_1}^2).
\end{align}
Consequently, we obtain
\begin{align*}
    0 & = \mathbb{E}_{x}^{\hat{V}_1}\left[\sum_{i=0}^{\hat\sigma_1 \wedge n-1}\varphi_1 (i) (g(X_i^1)-\lambda_1)+ \sum_{i=0}^{\hat\sigma_1 \wedge n-1}\varphi_1 (i) (w_1(X_{i+1}^1)-w_1(X_i^1))\right.\\
    & \phantom{=} \left.+ 1_{\{\hat\sigma_1\leq n\}} \varphi_1(\hat\sigma_1)c(X_{\hat\sigma_1}^1,X_{\hat\sigma_1}^2) + 1_{\{\hat\sigma_1\leq n\}} \varphi_1(\hat\sigma_1)(w_1(X_{\hat\sigma_1}^2)-  w_1(X_{\hat\sigma_1}^1) )\right].
\end{align*}
In fact, for any $k\in \mathbb{N}$, inductively we get
\begin{align*}
    0 & = \mathbb{E}_{x}^{\hat{V}_1}\left[\sum_{i=0}^{\hat\sigma_k \wedge n-1}\varphi_1 (i) (g(Y_i)-\lambda_1)+ \sum_{j=0}^{k-1}\sum_{i=\hat\sigma_{j}\wedge n}^{\hat\sigma_{j+1} \wedge n -1}\varphi_1 (i) (w_1(X_{i+1}^{j+1})-w_1(X_{i}^{j+1}))\right.\\
    & \phantom{,} \left.+ \sum_{j=1}^k 1_{\{\hat\sigma_j\leq n\}} \varphi_1(\hat\sigma_j)c(X_{\hat\sigma_j}^{j},X_{\hat\sigma_j}^{j+1}) + \sum_{j=1}^k 1_{\{\hat\sigma_j\leq n\}}\varphi_1(\hat\sigma_j) (w_1(X_{\hat\sigma_j}^{j+1})-w_1(X_{\hat\sigma_j}^j))\right].
\end{align*}
Next, using the fact that $\hat\sigma_k\to\infty$ as $k\to\infty$ (see the discussion following Equation~\eqref{eq:Fatou}), we obtain
\begin{align}\label{eq:th:joint_str_disc:1}
    0 & = \mathbb{E}_{x}^{\hat{V}_1}\left[\sum_{i=0}^{ n-1}\varphi_1 (i) (g(Y_i)-\lambda_1)+ \sum_{j=0}^\infty\sum_{i=\hat\sigma_{j}\wedge n}^{\hat\sigma_{j+1} \wedge n -1}\varphi_1 (i) (w_1(X_{i+1}^{j+1})-w_1(X_{i}^{j+1}))\right. \nonumber\\
    & \phantom{,} \left.+ \sum_{j=1}^\infty 1_{\{\hat\sigma_j\leq n\}} \varphi_1(\hat\sigma_j)c(X_{\hat\sigma_j}^{j},X_{\hat\sigma_j}^{j+1}) + \sum_{j=1}^\infty 1_{\{\hat\sigma_j\leq n\}}\varphi_1(\hat\sigma_j) (w_1(X_{\hat\sigma_j}^{j+1})-w_1(X_{\hat\sigma_j}^j))\right].
\end{align}
Next, let us define
\begin{multline}
    S:=\sum_{j=0}^\infty\sum_{i=\hat\sigma_{j}\wedge n}^{\hat\sigma_{j+1} \wedge n -1}\varphi_1 (i) (w_1(X_{i+1}^{j+1})-w_1(X_{i}^{j+1})) \\
    +\sum_{j=1}^\infty 1_{\{\hat\sigma_j\leq n\}}\varphi_1(\hat\sigma_j) (w_1(X_{\hat\sigma_j}^{j+1})-w_1(X_{\hat\sigma_j}^j))
\end{multline}
and note that
\begin{align*}
    S&= \sum_{j=0}^\infty \left(\sum_{i=\hat\sigma_{j}\wedge n}^{\hat\sigma_{j+1} \wedge n -2}w_1(X_{i+1}^{j+1})(\varphi_1 (i)-\varphi_1 (i+1))-1_{\{\hat\sigma_j\leq n\}} \varphi_1(\hat\sigma_{j})w_1(X_{\hat\sigma_j}^{j+1})\right.\nonumber\\
    &\phantom{=}\left. + 1_{\{\hat\sigma_j\leq n\}} \varphi_1(\hat\sigma_{j+1}\wedge n-1)w_1(X_{\hat\sigma_{j+1}\wedge n}^{j+1})\right)\\
    &\phantom{=}+ \sum_{j=1}^\infty 1_{\{\hat\sigma_j\leq n\}}\varphi_1(\hat\sigma_j) (w_1(X_{\hat\sigma_j}^{j+1})-w_1(X_{\hat\sigma_j}^j))\\
    &= -\varphi_1(0)w_1(X_{0}^{1})+\sum_{j=0}^\infty \sum_{i=\hat\sigma_{j}\wedge n}^{\hat\sigma_{j+1} \wedge n -2}w_1(X_{i+1}^{j+1})(\varphi_1 (i)-\varphi_1 (i+1))\nonumber\\
    &\phantom{=}+\sum_{j=0}^\infty 1_{\{\hat\sigma_{j+1}\leq n\}}w_1(X_{\hat\sigma_{j+1}}^{j+1})(\varphi_1(\hat\sigma_{j+1}-1)-\varphi_1(\hat\sigma_{j+1})) \\
    &\phantom{=}+\sum_{j=0}^\infty 1_{\{\hat\sigma_{j}\leq n<\hat\sigma_{j+1}\}}w_1(X_{n}^{j+1})\varphi_1(n-1).
\end{align*}
Consequently, we have
\begin{align*}
    S&  = -\varphi_1(0)w_1(X_{0}^{1})+\sum_{j=0}^\infty \sum_{i=\hat\sigma_{j}\wedge n}^{\hat\sigma_{j+1} \wedge n -1}w_1(X_{i+1}^{j+1})(\varphi_1 (i)-\varphi_1 (i+1))\nonumber\\
    &\phantom{=}+\sum_{j=0}^\infty 1_{\{\hat\sigma_{j}\leq n<\hat\sigma_{j+1}\}}w_1(X_{n}^{j+1})\varphi_1(n-1).
\end{align*}
Next, recalling the boundedness of $w_1$, we get
\begin{align}\label{eq:th:main_disc:0.45}
    |S|&\leq 2\Vert w_1\Vert \varphi_1(0)+\Vert w_1\Vert\sum_{j=0}^\infty \sum_{i=\hat\sigma_{j}\wedge n}^{\hat\sigma_{j+1} \wedge n -1}(\varphi_1 (i)-\varphi_1 (i+1))\nonumber\\
    & = 2\Vert w_1\Vert \varphi_1(0) + \Vert w_1\Vert\sum_{j=0}^\infty(\varphi_1 (\hat\sigma_{j}\wedge n)-\varphi_1 (\hat\sigma_{j+1}\wedge n))\nonumber\\
    & = 2\Vert w_1\Vert \varphi_1(0) + \Vert w_1\Vert\left(\varphi_1 (0)-\sum_{j=0}^\infty   1_{\{\hat\sigma_j\leq n<\hat\sigma_{j+1}\}} \varphi_1(n)\right)\leq 3\Vert w_1\Vert \varphi_1(0).
\end{align}
Thus, as in~\eqref{eq:pr:ver_th_disc_inc:2.5}, dividing~\eqref{eq:th:joint_str_disc:1} by $\sum_{i=0}^{n-1}\varphi_1(i)$ and letting $n\to\infty$, we have
\begin{align}\label{eq:th:main_disc:0.5}
    \lambda_1 & = \limsup_{n\to\infty}\frac{\mathbb{E}_{x}^{\hat{V}_1}\left[\sum_{i=0}^{n-1}\varphi_1 (i) g(Y_i)+ \sum_{i=1}^\infty 1_{\{\hat\sigma_i\leq n\}} \varphi_1(\hat\sigma_i)c(Y_{\hat\sigma_i^-},Y_{\hat\sigma_i}) \right]}{\sum_{i=0}^{n-1}\varphi_1(i)} \nonumber\\
    & = J_1(x,\hat{V}).
\end{align}

Second, let $V:=(\tau_i,\xi_i)_{i=1}^\infty\in \mathbb{V}_1$ and note that using~\eqref{eq:Bellman_disc_inc}, for any $x\in E$, $n\in \mathbb{N}$, and $i=0, 1, \ldots, \tau_1 \wedge n-1$, we obtain
\begin{equation}\label{eq:th_main_disc:3}
  \varphi_1 (i) w_1(X_i^1)\leq \varphi_1 (i) (g(X_i^1)-\lambda_1)+\varphi_1 (i) \mathbb{E}_{X_i^1}^{V}[w_1(X_1)], \quad \mathbb{P}_{x}^{V} \, \text{a.s.}  
\end{equation}
Also, on the event $\{\tau_1<\infty\}$, we have
\begin{align}\label{eq:th_main_disc:4}
\varphi_1(\tau_1)w_1(X_{\tau_1}^1) &\leq \varphi_1(\tau_1)Mw_1(X_{\tau_1}^1) \nonumber\\
&\leq \varphi_1(\tau_1)c(X_{\tau_1}^1,\xi_1)+\varphi_1(\tau_1)w_1(\xi_1).
\end{align}
Thus, as in~\eqref{eq:th:joint_str_disc:1}, we get
\begin{align}\label{eq:th:main_disc:0.55}
    0 & \leq \mathbb{E}_{x}^{V}\left[\sum_{i=0}^{ n-1}\varphi_1 (i) (g(Y_i)-\lambda_1)+ \sum_{j=0}^\infty\sum_{i=\tau_{j}\wedge n}^{\tau_{j+1} \wedge n -1}\varphi_1 (i) (w_1(X_{i+1}^{j+1})-w_1(X_{i}^{j+1}))\right. \nonumber\\
    & \phantom{,} \left.+ \sum_{j=1}^\infty 1_{\{\tau_{j}\leq n\}} \varphi_1(\tau_{j})c(X_{\tau_{j}}^{j},X_{\tau_{j}}^{j+1}) + \sum_{j=1}^\infty 1_{\{\tau_{j}\leq n\}}\varphi_1(\tau_{j}) (w_1(X_{\tau_{j}}^{j+1})-w_1(X_{\tau_{j}}^j))\right].
\end{align}
Next, using the argument leading to~\eqref{eq:th:main_disc:0.45}, we get
\begin{multline*}
    \left|\sum_{j=0}^\infty\sum_{i=\tau_{j}\wedge n}^{\tau_{j+1} \wedge n -1}\varphi_1 (i) (w_1(X_{i+1}^{j+1})-w_1(X_{i}^{j+1}))\right.\\
    \left.+\sum_{j=1}^\infty 1_{\{\tau_{j}\leq n\}}\varphi_1(\tau_{j}) (w_1(X_{\tau_{j}}^{j+1})-w_1(X_{\tau_{j}}^j))\right|\leq 3\Vert w_1\Vert \varphi(0).
\end{multline*}
Thus, dividing~\eqref{eq:th:main_disc:0.55} by $\sum_{i=0}^{n-1}\varphi_1(i)$ and letting $n\to\infty$, we have
\[
\lambda_1\leq J_1^d(x,V). 
\]
Combining this with~\eqref{eq:th:main_disc:0.5}, we obtain
\[
 \inf_{V\in \mathbb{V}_1} J^d_1(x,V) =J^d_1(x,\hat{V}_1) =\lambda_1, \quad x\in E.
\]
Recalling Corollary~\ref{cor:opt_str_disc_undisc}, we conclude the proof.
\end{proof}

In Theorem~\ref{th:main_discrete} we showed that the values of the functionals $J_h^d$ and $J_h$ are the same for the optimal strategy $\hat{V}_h$ given by~\eqref{eq:opt_str_disc_cons}. Now, we extend this result to the family of any Markovian stationary strategies. More specifically, let $D\in \mathcal{E}$ and let $\psi\colon D^c\to U$ be a bounded measurable map satisfying $\psi(D^c)\subset D$. Then, we define a strategy $\bar{V}:=(\bar\sigma_i,\bar\xi_i)_{i=1}^\infty$ by
\begin{align}\label{eq:str_Markov}
    \bar\sigma_i&:=h\inf\{n\in \mathbb{N}\colon n\geq \bar\sigma_{i-1}/h, \, X_{nh}^i\notin D\}, \nonumber\\
    \bar\xi_i&:=\psi(X_{\bar\sigma_i}^i)1_{\{\bar\sigma_i<\infty\}}+\xi_0 1_{\{\bar\sigma_i=\infty\}},
\end{align}
where $\bar\sigma_0\equiv 0$ and $\xi_0\in U$ is some fixed point. Note that this is a strategy with impulse times defined as the consecutive exit times for the set $D$ and after-impulse states given by the map $\psi$. Also, the condition $\psi(D^c)\subset D$ prevents the situation when applying an impulse triggers another immediate impulse.

With the strategy $\bar{V}$, we associate the Poisson equation of the form
\begin{equation}\label{eq:Poisson}
\begin{split}
    w_{\bar{V}}(x)&=\mathbb{E}_x\left[\sum_{i=0}^{(\tau_{D}/h)\wedge n-1}(g(X_{ih})-\lambda_{\bar{V}})+w_{\bar{V}}(X_{(nh)\wedge \tau_D})\right], \quad x\in E,\\
    w_{\bar{V}}(x)&=\frac{1}{h}c(x,\psi(x))+w_{\bar{V}}(\psi(x)), \quad x\notin D,
\end{split}
\end{equation}
where $\tau_D:=h\inf\{n\in \mathbb{N}\colon X_{nh}\notin D\}$. Note that this equation could be seen as a version of~\eqref{eq:Bellman_disc_inc} with a fixed control. In the following proposition we show that there exists a solution $(w_{\bar{V}}, \lambda_{\bar{V}})$ to~\eqref{eq:Poisson}.
\begin{proposition}\label{pr:Poisson}
    There exists a bounded measurable map $w_{\bar{V}}\colon E\to \mathbb{R}$ and a constant $\lambda_{\bar{V}}$ such that the pair $(w_{\bar{V}}, \lambda_{\bar{V}})$ satisfies~\eqref{eq:Poisson}.
\end{proposition}
\begin{proof}
The argument is partially based on the proof of Proposition~\ref{pr:existence_discounted}, thus we provide only an outline. First, using the span-contraction approach, we find a bounded measurable map $w_{\bar{V}}\colon D\to \mathbb{R}$ and a constant $\lambda_{\bar{V}}\in \mathbb{R}$ satisfying
\begin{equation}\label{eq:pr:Poisson:1}
    w_{\bar{V}}(x)=g(x)-\lambda_{\bar{V}}+\mathbb{E}_x[w_{\bar{V}}(X_h)], \quad x\in D.
\end{equation}
More specifically, we consider the family $\mathcal{B}_b(D)$ of bounded measurable real-valued functions on $D$ with the span semi-norm $\Vert f\Vert_{sp}^D:=\sup_{x,y\in D}(f(x)-f(y))$, $f\in \mathcal{B}_b(D)$. Also, for any $v\in \mathcal{B}_b(D)$ and $x\in D$, we define the operator
\[
Fv(x):=g(x)+\mathbb{E}_x[1_{\{X_h\in D\}}v(X_h)+1_{\{X_h\in D^c\}}(c(X_h,\psi(X_h))+v(\psi(X_h)))];
\]
note that directly from the definition we get $F\mathcal{B}_b(D)\subset \mathcal{B}_b(D)$. Next, recalling that $\psi(D^{c})\subset D$ and using Assumption~\eqref{A2}, as in~\eqref{eq:pr:existence_discounted:span}, for any $v_1,v_2\in \mathcal{B}_b(D)$, we obtain
\begin{align*}
    \Vert Fv_1-Fv_2\Vert_{sp}^D  & = \sup_{x,y\in D} \left(\int_D v_1(z)P_h(x,dz) + \int_{D^c}\frac{1}{h}(c(z,\psi(z))+v_1(\psi(z)))P_1(x,dz) \right.\nonumber\\
    & \phantom{=} -  \int_D v_2(z)P_h(x,dz) - \int_{D^c}(\frac{1}{h}c(z,\psi(z))+v_2(\psi(z)))P_1(x,dz)\nonumber\\
    & \phantom{=} - \int_D v_1(z)P_h(y,dz) - \int_{D^c}(\frac{1}{h}c(z,\psi(z))-v_1(\psi(z)))P_1(y,dz)\nonumber\\
    & \phantom{=} \left.+ \int_D v_2(z)P_h(y,dz) + \int_{D^c}(\frac{1}{h}c(z,\psi(z))+v_2(\psi(z)))P_1(y,dz)\right)\nonumber\\
    &=\sup_{x,y\in D}\left(\int_D (v_1(z)-v_2(z))P_h(x,dz)-\int_D (v_1(z)-v_2(z))P_h(y,dz)\right.\nonumber\\
    &\phantom{=}+\int_{D^c}(v_1(\psi(z))-v_2(\psi(z))) P_h(x,dz)\nonumber\\
    &\phantom{=}\left.-\int_{D^c}(v_1(\psi(z))-v_2(\psi(z))) P_h(y,dz)\right)\nonumber\\
    &\leq \Lambda_1\Vert v_1-v_2\Vert_{sp}^D .
\end{align*}
Thus, from~\eqref{A2}, we conclude that the operator $F$ is a contraction with respect to the span semi-norm on $\mathcal{B}_b(D)$. Hence, using the Banach fixed-point theorem, we get that there exists a map $w_{\bar{V}}\in \mathcal{B}_{b}(D)$ and a constant $\lambda_{\bar{V}}$ such that $Fw_{\bar{V}} (x)= \lambda_{\bar{V}}+w_{\bar{V}}(x)$, $x\in D$. 

Next, we extend the definition of $w_{\bar{V}}$ to the full space $E$ by setting
\begin{equation}\label{eq:pr:Poisson:2}
    w_{\bar{V}}(x):=\frac{1}{h}c(x,\psi(x))+w_{\bar{V}}(\psi(x)), \quad x\notin D;
\end{equation}
note that the definition is correct since $\psi(D^c)\subset D$, and $w_{\bar{V}}$ at the right hand-side is already defined there. Consequently, recalling the definition of $F$ and~\eqref{eq:pr:Poisson:2}, we get that~\eqref{eq:pr:Poisson:1} is satisfied. 

Finally, combining~\eqref{eq:pr:Poisson:1} and~\eqref{eq:pr:Poisson:2}, for any $x\in E$ and $n\in \mathbb{N}$, we get
\begin{align*}
    w_{\bar{V}}(x)&=\mathbb{E}_x\Bigg[\sum_{i=0}^{(\tau_{D}/h)\wedge n-1}(g(X_{ih})-\lambda_{\bar{V}})\\
    &\phantom{===}+1_{\{\tau_D\leq nh\}}(\frac{1}{h}c(X_{\tau_D},\psi(X_{\tau_D}))+w_{\bar{V}}(\psi(X_{\tau_D})))+1_{\{\tau_D> nh\}}w_{\bar{V}}(X_{nh})\Bigg]\\
    &=\mathbb{E}_x\left[\sum_{i=0}^{(\tau_{D}/h)\wedge n-1}(g(X_{ih})-\lambda_{\bar{V}})+w_{\bar{V}}(X_{(nh)\wedge \tau_D})\right],
\end{align*}
which concludes the proof.
\end{proof}

Using Proposition~\ref{pr:Poisson}, we show that the values of the functionals $J^d_h$ and $J_h$ are equal for  Markovian stationary strategies.

\begin{theorem}\label{th:same_payoff_disc}
    Let the strategy $\bar{V}$ be given by~\eqref{eq:str_Markov} and let $(w_{\bar{V}}, \lambda_{\bar{V}})$ be a solution to~\eqref{eq:Poisson}. Also, let $J^d_h$ and $J_h$ be given by~\eqref{eq:funct_time_inc_disc} and~\eqref{eq:funct_time_cons_disc}, respectively. Then, we get
    \[
J^d_h(x,\bar{V})=\lambda_{\bar{V}}=J_h(x,\bar{V}), \quad x\in E.
    \]
\end{theorem}
\begin{proof} The argument is based on the proofs of Proposition~\ref{pr:opt_str_disc_disc} and Theorem~\ref{th:main_discrete}, thus we provide only an outline.
First, note that as in Proposition~\ref{pr:martingale_disc_inc}, we may show that the process
\[
z_{\bar{V}}(n):=\sum_{i=0}^{(\tau_{D}/h)\wedge n-1}(g(X_{ih})-\lambda_{\bar{V}})+w_{\bar{V}}(X_{(nh)\wedge \tau_D}), \quad n\in \mathbb{N},
\]
is a martingale. Thus, as in the first part of the proof of Proposition~\ref{pr:opt_str_disc_disc}, we show $\lambda_{\bar{V}}=J_h(x,\bar{V})$, $x\in E$. 

Second, note that setting $n=1$ in~\eqref{eq:Poisson}, we obtain
\[
w_{\bar{V}}(x)=g(x)-\lambda_{\bar{V}}+\mathbb{E}_x[w_{\bar{V}}(X_h)], \quad x\in D.
\]
Thus, for any $x\in E$, $n\in \mathbb{N}$, and $i=0, 1, \ldots, (\bar\sigma_1/h) \wedge n-1$, we obtain
\begin{equation}\label{eq:th:same_payoff_disc:1}
  \varphi_1 (i) w_{\bar{V}}(X_i^1)= \varphi_1 (i) (g(X_i^1)-\lambda_{\bar{V}})+\varphi_1 (i) \mathbb{E}_{X_i^1}^{\bar{V}}[w_{\bar{V}}(X_1)], \quad \mathbb{P}_{x}^{\bar{V}} \, \text{a.s.} 
\end{equation}
Also, directly from~\eqref{eq:Poisson}, on the event $\{\bar\sigma_1/h<\infty\}$, we have
\begin{align*}
\varphi_1(\bar\sigma_1/h)w_{\bar{V}}(X_{\bar\sigma_1}^1) = \varphi_1(\bar\sigma_1/h)c(X_{\hat\sigma_1}^1,\bar\xi_1)+\varphi_1(\bar\sigma_1/h)w_{\bar{V}}(\bar\xi_1).
\end{align*}
Comparing these two identities with~\eqref{eq:th_main_disc:1} and~\eqref{eq:th_main_disc:2}, and using the argument from the first part of the proof of Theorem~\ref{th:main_discrete}, we obtain $J^d_h(x,\bar{V})=\lambda_{\bar{V}}$, $x\in E$, which concludes the proof.
\end{proof}


\section{Equivalence of discounted and undiscounted problems}\label{S:continuous}

In this section we show that the optimal values of continuous time discounted and undiscounted impulse control problems are equal. This could be seen as an extension of Theorem~\ref{th:main_discrete} stated for the discrete-time setting. Note that, in the continuous time framework, we do not have a simple one-step-ahead Bellman equation~\eqref{eq:Bellman_disc_cons} and some additional arguments are needed. In fact, in the proof we use two types of approximation schemes for the continuous time problems: one with discrete strategies and one with discrete-type cost functionals. Note that this could also be used to numerically approximate the problems. Also, in this way we do not need to solve the associated continuous time Bellman equation; see e.g.~\cite{PalSte2017} and references therein for an exemplary class of assumptions used to solve it.

First, we show an alternative representation of~\eqref{eq:funct_time_inc}. This result is extensively used in the subsequent analyses.

\begin{lemma}\label{lm:alt_form}
Let $J^d$ be given by~\eqref{eq:funct_time_inc}. Then, for any $h>0$, $x\in E$, $V\in \mathbb{V}$, we get
\begin{equation}\label{eq:funct_time_inc:alt}
    J^d(x,V)=\limsup_{n\to\infty}\frac{\mathbb{E}_{x}^V\left[\int_0^{nh} \beta(s)g(Y_{s})ds+\sum_{i=1}^\infty 1_{\{\tau_i\leq nh\}}\beta(\tau_i)c(X_{\tau_i}^i,X_{\tau_i}^{i+1}) \right]}{\int_0^{nh}\beta(s)ds}.
\end{equation}
\end{lemma}
\begin{proof}
    Denote by $\hat{J}^d(x,V)$ the right hand side of~\eqref{eq:funct_time_inc:alt} and note that by considering the sequence $(nh)_{n\in \mathbb{N}}$, we get $J^d(x,V)\geq \hat{J}^d(x,V)$. For the reverse inequality, note that using the boundedness of $g$ and the non-negativity of $c$, for any $T>0$, $h>0$, $x\in E$, and $V\in \mathbb{V}$, we get
    \begin{multline*}
     \mathbb{E}_{x}^V\left[\int_0^{T} \beta(s)g(Y_{s})ds+\sum_{i=1}^\infty 1_{\{\tau_i\leq T\}}\beta(\tau_i)c(X_{\tau_i}^i,X_{\tau_i}^{i+1}) \right]\\
    \leq \mathbb{E}_{x}^V\left[\int_0^{([T/h]+1)h} \beta(s)g(Y_{s})ds+\sum_{i=1}^\infty 1_{\{\tau_i\leq ([T/h]+1)h\}}\beta(\tau_i)c(X_{\tau_i}^i,X_{\tau_i}^{i+1}) \right] +h\Vert g\Vert, 
    \end{multline*}
where $[a]:=\sup\{k\in \mathbb{Z}\colon k\leq a\}$ denotes the integer part of $a\in \mathbb{R}$. Thus, diving both hand sides by $\int_0^T \beta(s)ds$ and noting that $\lim_{T\to\infty}\frac{\int_0^{([T/h]+1)h} \beta(s)ds}{\int_0^T \beta(s)ds}=1$, we get
\begin{multline*}
    J^d(x,V)\\
    \leq \limsup_{T\to\infty} \frac{\mathbb{E}_{x}^V \left[\int_0^{([T/h]+1)h} \beta(s)g(Y_{s})ds+\sum_{i=1}^\infty 1_{\{\tau_i\leq ([T/h]+1)h\}}\beta(\tau_i)c(X_{\tau_i}^i,X_{\tau_i}^{i+1}) \right]}{\int_0^{([T/h]+1)h}\beta(s)ds}\\
    \leq \hat{J}^d(x,V),
\end{multline*}
which concludes the proof.
\end{proof}

Next, we show that the continuous time control problem with generalised discounting and discrete control strategies could be approximated by the problem with discrete strategies and discrete-type optimality functional.

\begin{theorem}\label{th:approx_disc}
Let $J^d$ and $J^d_h$ be given by~\eqref{eq:funct_time_inc} and~\eqref{eq:funct_time_inc_disc}, respectively. Then, we get 
    \[\lim_{h\to 0}|\inf_{V\in \mathbb{V}_h} J^d_h(x,V)-\inf_{V\in \mathbb{V}_h} J^d(x,V)|=0.
    \]
\end{theorem}
\begin{proof}
First, note that, using~\eqref{eq:phi_h} and~\eqref{eq:funct_time_inc_disc}, we get that $J^d_h$ satisfies
\begin{align}\label{eq:th:approx_disc:funct_time_inc_disc}
    J^d_h(x,V)&=\limsup_{n\to\infty}\frac{1}{\frac{1}{h}\int_{0}^{nh}\beta(s)ds}\mathbb{E}_{x}^V\left[\sum_{i=0}^{n-1} \frac{1}{h}\int_{ih}^{(i+1)h}\beta(s)ds\, g(Y_{ih})\right.\nonumber\\
    &\phantom{=}\left.+\sum_{i=1}^\infty 1_{\{\tau_i\leq nh\}}\frac{1}{h}\int_{\tau_i}^{\tau_i+h}\beta(s)ds \frac{1}{h}c(X_{\tau_i}^i,X_{\tau_i}^{i+1}) \right]\nonumber\\
    &=\limsup_{n\to\infty}\frac{1}{\int_{0}^{nh}\beta(s)ds}\mathbb{E}_{x}^V\left[\sum_{i=0}^{n-1} \int_{ih}^{(i+1)h}\beta(s)ds \,g(Y_{ih})\right.\nonumber\\
    &\phantom{=}\left.+\sum_{i=1}^\infty 1_{\{\tau_i\leq nh\}}\frac{1}{h}\int_{\tau_i}^{\tau_i+h}\beta(s)ds\, c(X_{\tau_i}^i,X_{\tau_i}^{i+1}) \right].
\end{align}
Next, by considering the \textit{no-impulse} strategy, we obtain
\[
\inf_{V\in \mathbb{V}_h}J^d(x,V)\leq \Vert g\Vert \quad \text{and} \quad \inf_{V\in \mathbb{V}_h}J^d_h(x,V)\leq \Vert g\Vert.
\]
Now, let us denote by $\mathbb{V}_h^* \subset \mathbb{V}_h$ the subfamily of strategies $V$ for which we have $J^d(x,V)\leq \Vert g\Vert$ or $J^d_h(x,V)\leq \Vert g\Vert$. Thus, for any $h>0$ and $x\in E$, we get
\begin{align*}
    |\inf_{V\in \mathbb{V}_h} J^d(x,V)-\inf_{V\in \mathbb{V}_h} J^d_h(x,V)|&=    |\inf_{V\in \mathbb{V}_h^*} J^d(x,V)-\inf_{V\in \mathbb{V}_h^*} J^d_h(x,V)|\\
    &\leq \sup_{V\in \mathbb{V}_h^*}|J^d(x,V)-J^d_h(x,V)|.
\end{align*}
Hence, for any  $h>0$, $x\in E$, and $V\in \mathbb{V}_h^*$, using Lemma~\ref{lm:alt_form}, we have
\begin{multline}\label{eq:th:approx_disc:1}
    |J^d(x,V)-J^d_h(x,V)|\leq \\
    \limsup_{n\to\infty}\frac{1}{\int_0^{nh}\beta(s)ds}\left| \mathbb{E}_{x}^V\left[\int_0^{nh} \beta(s)g(Y_{s})ds -
    \sum_{i=0}^{n-1} \int_{ih}^{(i+1)h}\beta(s)ds \,g(Y_{ih})\right.\right.\\
    \left.\left. +\sum_{i=1}^\infty 1_{\{\tau_i\leq nh\}}c(X_{\tau_i}^i,X_{\tau_i}^{i+1})\left(\beta(\tau_i)-\frac{1}{h}\int_{\tau_i}^{\tau_i+h}\beta(s)ds  \right)\right]\right|.
\end{multline}
Next, note that using the compactness of $E$, Assumption~\eqref{Feller}, and Theorem 3.7 from~\cite{Dyn1965}, we have
\begin{equation*}
   \sup_{y\in E}\sup_{t\in [0,h]} |\mathbb{E}_y[g(X_t)]-g(y)|=:\delta(h) \to 0, \quad h\downarrow 0; 
\end{equation*}
see also theorem on page 25 in~\cite{Mey1967} and Proposition 2.4 in Chapter III in~\cite{RevYor1999}. Thus, recalling the facts that the controlled process is Markov (with the uncontrolled dynamics) between the consecutive impulses and the shifts can be applied only at the multiplicities of $h$, we obtain
\begin{multline}\label{eq:th:approx_disc:1.5}
    \left|\mathbb{E}_{x}^V\left[\int_0^{nh} \beta(s)g(Y_{s})ds -
    \sum_{i=0}^{n-1} \int_{ih}^{(i+1)h}\beta(s)ds \,g(Y_{ih})\right]\right|\\
    = \left|\sum_{i=0}^{n-1} \mathbb{E}_{x}^V\left[\int_{ih}^{(i+1)h}\beta(s)(g(Y_s)-g(Y_{ih}))ds \right]\right|\\
    \leq \delta(h)\sum_{i=0}^{n-1} \int_{ih}^{(i+1)h}\beta(s)ds  = \delta(h) \int_0^{nh}\beta(s)ds.
\end{multline}
Next, we show that, for $h>0$ small enough and any $V\in \mathbb{V}_h^*$, we obtain 
\begin{equation}\label{eq:th:approx_disc:2}
    \limsup_{n\to\infty}\frac{\mathbb{E}_{x}^V\left[ \sum_{i=1}^\infty 1_{\{\tau_i\leq nh\}}\beta(\tau_i)\right]}{\int_0^{nh}\beta(s)ds}\leq \frac{3\Vert g\Vert}{c_0}.
\end{equation}
Indeed, if $J^d_h(x,V)\leq \Vert g\Vert$, then we have
\begin{align*}
    \Vert g\Vert \geq J^d_h(x,V)\geq \limsup_{n\to\infty}\frac{-\Vert g\Vert\int_0^{nh} \beta(s)ds+\mathbb{E}_x^V[\sum_{i=1}^\infty 1_{\{\tau_i\leq nh\}}c_0\beta(\tau_i+h)]}{\int_0^{nh}\beta(s)ds}.
\end{align*}
Thus, using~\eqref{eq:beta_mult} to get $\beta(\tau_i+h)\geq \beta(\tau_i)\beta(h)$, for $h$ small enough, we obtain
\[
\limsup_{n\to\infty}\frac{\mathbb{E}_{x}^V\left[ \sum_{i=1}^\infty 1_{\{\tau_i\leq nh\}}\beta(\tau_i)\right]}{\int_0^{nh}\beta(s)ds}\leq \frac{2\Vert g\Vert}{c_0\beta(h)}\leq \frac{3\Vert g\Vert}{c_0}.
\]
If $J^d(x,V)\leq \Vert g\Vert$, then, using a similar argument, one may show
\[
\limsup_{n\to\infty}\frac{\mathbb{E}_{x}^V\left[ \sum_{i=1}^\infty 1_{\{\tau_i\leq nh\}}\beta(\tau_i)\right]}{\int_0^{nh}\beta(s)ds}\leq \frac{2\Vert g\Vert}{c_0}\leq \frac{3\Vert g\Vert}{c_0},
\]
which concludes the proof of~\eqref{eq:th:approx_disc:2}.

Now, using the monotonicity of $\beta$ and~\eqref{eq:beta_mult}, for any $t,h\geq 0$, we get
\begin{align}
    0\leq \beta(t)-\frac{1}{h}\int_{t}^{t+h}\beta(s)ds\leq \beta(t)-\beta(t+h)\leq \beta(t)(1-\beta(h)).
\end{align}
Thus, combining this with~\eqref{eq:th:approx_disc:1},~\eqref{eq:th:approx_disc:1.5}, and~\eqref{eq:th:approx_disc:2}, we have
\begin{multline*}
    |J^d(x,V)-J^d_h(x,V)|\leq \\
    \limsup_{n\to\infty}\frac{\delta(h)\int_0^{nh} \beta(s)ds +(1-\beta(h))\mathbb{E}_{x}^V\left[ \sum_{i=1}^\infty \beta(\tau_i)1_{\{\tau_i\leq nh\}}c(X_{\tau_i}^i,X_{\tau_i}^{i+1})\right]}{\int_0^{nh}\beta(s)ds}\\
    \leq \delta(h)+(1-\beta(h))\Vert c\Vert \frac{3\Vert g\Vert}{c_0} .
\end{multline*}
Hence, letting $h\to 0$, we conclude the proof.
\end{proof}

The next result could be seen as a version of Theorem~\ref{th:approx_disc} for the undiscounted functional.

\begin{theorem}\label{th:approx_undisc}
Let $J$ and $J_h$ be given by~\eqref{eq:funct_time_cons} and~\eqref{eq:funct_time_cons_disc}, respectively. Then, we get 
    \[\lim_{h\to 0}|\inf_{V\in \mathbb{V}_h} J_h(x,V)-\inf_{V\in \mathbb{V}_h} J(x,V)|=0.
    \]
\end{theorem}
\begin{proof}
Note that, for any $h>0$,  $x\in E$, and $V\in \mathbb{V}$, as in Lemma~\ref{lm:alt_form}, we may show
\begin{align*}
    J(x,V)&=\limsup_{n\to\infty}\frac{1}{nh}\mathbb{E}_{x}^V\left[\int_0^{nh} g(Y_s)ds+\sum_{i=1}^\infty 1_{\{\tau_i\leq nh\}}c(X_{\tau_i}^i,X_{\tau_i}^{i+1}) \right].
\end{align*}
Also, directly from~\eqref{eq:funct_time_cons_disc}, we obtain
\[
    J_h(x,V)=\limsup_{n\to\infty}\frac{1}{nh}\mathbb{E}_{x}^V\left[\sum_{i=1}^n hg(Y_{ih})+\sum_{i=1}^\infty 1_{\{\tau_i\leq nh\}}c(X_{\tau_i}^i,X_{\tau_i}^{i+1}) \right].
\]
Thus, using the argument exploited in Theorem~\ref{th:approx_disc} with $\beta\equiv 1$, we may conclude the proof. For brevity, we omit the details.
\end{proof}

Now, we show that the problem with generalised discounting could be approximated by the problem with discrete strategies, where the impulse times are restricted to the time grid spanned by $h>0$.

\begin{theorem}\label{th:approx_dyadic_disc}
Let $J^d$ be given by~\eqref{eq:funct_time_inc}. Then, we get 
    \[\lim_{h\to 0}\inf_{V\in \mathbb{V}_h} J^d(x,V)=\inf_{V\in \mathbb{V}} J^d(x,V).
    \]
\end{theorem}
\begin{proof} 
First, note that, for any $h>0$ and $x\in E$, we have $\inf_{V\in \mathbb{V}_h} J^d(x,V)\geq\inf_{V\in \mathbb{V}} J^d(x,V)$. We show that in the limit we have equality. To see this, let $x\in E$, $\varepsilon>0$, and $V^{\varepsilon}:=(\tau_i,\xi_i)_{i=1}^\infty\in \mathbb{V}$ be an $\varepsilon$-optimal strategy for $\inf_{V\in \mathbb{V}} J^d(x,V)$. Also, for any $h>0$, let us define a strategy $V^{\varepsilon}_h:=(\tau_i^h,\xi_i^h)_{i=1}^\infty$ as
\begin{align}\label{eq:th:approx_dyadic_disc:strategy}
\tau_i^h:=h\inf\{k\in \mathbb{N}\colon \tau_i\leq kh\}, \quad \xi_i^h&:=Y_{\tau_i^h}, \quad i=1,2,\ldots
\end{align}
Note that the strategy $V^{\varepsilon}_h$ is simply a $\mathbb{V}_h$-approximation of $V^{\varepsilon}$ with the impulse times in $\{0,h, 2h, \ldots\}$ and the impulses to the current state of the process controlled by $V^{\varepsilon}$. In the following, to ease the notation, by $Y$ and $Z$, we denote the processes controlled by $V^{\varepsilon}$ and $V^{\varepsilon}_h$, respectively. Recall that we have $Y_t:=X_t^i$ for $t\in [\tau_{i-1},\tau_i)$. Using the same coordinate process $(X_t^1, X_t^2, \ldots)$, we construct the process $Z$ as $Z_t:=X^i_{t}$, $t\in [\tau_{i-1}^h,\tau_i^h)$. Finally, without any loss of generality, we assume that $g(\cdot)\geq 0$ as for generic $g\in \mathcal{C}_b(E)$ we may set $\tilde{g}(\cdot):=g(\cdot)-\Vert g\Vert$ and the $J^d(x,V)$ changes by a constant.

Within this framework, using Lemma~\ref{lm:alt_form}, we have
\begin{align}\label{eq:th:approx_dyadic_disc:1}
    0&\leq \inf_{V\in \mathbb{V}_h} J^d(x,V)-\inf_{V\in \mathbb{V}} J^d(x,V)\leq J^d(x,V^{\varepsilon}_h)-J^d(x,V^{\varepsilon})+\varepsilon \nonumber\\
    & \leq  \limsup_{n\to\infty}\frac{1}{\int_0^{nh}\beta(s)ds}\left|\mathbb{E}_{x}^{V^{\varepsilon}}\left[\int_0^{nh} \beta(s)g(Y_{s})ds+\sum_{i=1}^\infty 1_{\{\tau_i\leq nh\}}\beta(\tau_i)c(X_{\tau_i}^i,X_{\tau_i}^{i+1}) \right]\right.\nonumber\\
    &\left.\phantom{=}-\mathbb{E}_{x}^{V^{\varepsilon}}\left[\int_0^{nh} \beta(s)g(Z_{s})ds+\sum_{i=1}^\infty 1_{\{\tau_i^h\leq nh\}}\beta(\tau_i^h)c(X_{\tau_{i}^h}^i,X_{\tau_i^h}^{i+1}) \right]\right|+\varepsilon.
\end{align}
Next, note that $Y_t=Z_t$ for $t\in [\tau_{i-1}^h,\tau_{i})$ and $\tau_{i}^h-\tau_{i}\leq h$. 
Thus, we obtain
\begin{multline}\label{eq:th:approx_dyadic_disc:1.2}
    \left|\mathbb{E}_{x}^{V^{\varepsilon}}\left[\int_0^{nh} \beta(s)g(Y_{s})ds\right]-\mathbb{E}_{x}^{V^{\varepsilon}}\left[\int_0^{nh} \beta(s)g(Z_{s})ds\right]\right|\\
    \leq  2\Vert g\Vert \mathbb{E}^{V^{\varepsilon}}_x\left[\sum_{i=1}^{\infty}   1_{\{\tau_i\leq nh\}}\int_{\tau_i}^{\tau_i^h}\beta(s)ds\right] \\
    \leq 2\Vert g\Vert h\mathbb{E}^{V^{\varepsilon}}_x\left[\sum_{i=1}^{\infty}   1_{\{\tau_i\leq nh\}}\beta(\tau_i)\right].
\end{multline}
Next, noting that $1_{\{\tau_i\leq nh\}}=1_{\{\tau_i^h\leq nh\}}$, we get
\begin{align}\label{eq:th:approx_dyadic_disc:1.5}
   \bigg|\mathbb{E}_{x}^{V^{\varepsilon}}&\left[\sum_{i=1}^\infty 1_{\{\tau_i\leq nh\}}\beta(\tau_i)c(X_{\tau_i}^i,X_{\tau_i}^{i+1}) \right]- \mathbb{E}_{x}^{V^{\varepsilon}}\left[\sum_{i=1}^\infty 1_{\{\tau_i^h\leq nh\}}\beta(\tau_i^h)c(X_{\tau_{i}^h}^i,X_{\tau_i^h}^{i+1})  \right]\bigg| \nonumber\\
   &\leq \mathbb{E}_{x}^{V^{\varepsilon}}\left[\sum_{i=1}^\infty 1_{\{\tau_i\leq nh\}}\beta(\tau_i)\left|c(X_{\tau_i}^i,X_{\tau_i}^{i+1})-c(X_{\tau_{i}^h}^i,X_{\tau_i}^{i+1})  \right| \right] \nonumber\\  
   &\phantom{=} + \mathbb{E}_{x}^{V^{\varepsilon}}\left[\sum_{i=1}^\infty 1_{\{\tau_i\leq nh\}}\beta(\tau_i^h)\left|c(X_{\tau_i^h}^i,X_{\tau_i^h}^{i+1})-c(X_{\tau_{i}^h}^i,X_{\tau_i}^{i+1})  \right| \right] \nonumber\\
   &\phantom{=} + \mathbb{E}_{x}^{V^{\varepsilon}}\left[\sum_{i=1}^\infty 1_{\{\tau_i\leq nh\}}(\beta(\tau_i)-\beta(\tau_i^h))c(X_{\tau_{i}^h}^i,X_{\tau_i}^{i+1})   \right].
\end{align}
Also, using the compactness of $E$ and the continuity of $c$, we may find a continuous and bounded map $K\colon \mathbb{R}\to\mathbb{R}$ such that $K(0)=0$ and 
\begin{align*}
    \sup_{\xi\in E}|c(x_1,\xi)-c(x_2,\xi)|&\leq K(\rho(x_1,x_2)), \quad x_1,x_2\in E,\\
    \sup_{y\in E}|c(y,\xi_1)-c(y,\xi_2)|&\leq K(\rho(\xi_1,\xi_2)), \quad \xi_1,\xi_2\in E.
\end{align*}
Let $r>0$ be small enough to get $K(z)\leq h$ for $|z|<r$. Using the compactness of $E$, Assumption~\eqref{Feller}, and Proposition 6.4 in~\cite{BasSte2018}, we have
\begin{equation*}
   \sup_{y\in E} \mathbb{P}_y[\sup_{t\in [0,h]}\rho(X_t,X_0)\geq r]=:\gamma(h) \to 0, \quad h\downarrow 0. 
\end{equation*}
Then, recalling that the controlled process between impulses is Markov with the original (uncontrolled) dynamics, we obtain
\begin{align*}
    \mathbb{E}_{x}^{V^{\varepsilon}}&\left[\sum_{i=1}^\infty 1_{\{\tau_i\leq nh\}}\beta(\tau_i)\left|c(X_{\tau_i}^i,X_{\tau_i}^{i+1})-c(X_{\tau_{i}^h}^i,X_{\tau_i}^{i+1})  \right| \right]\\
    & \leq \mathbb{E}_{x}^{V^{\varepsilon}}\left[\sum_{i=1}^\infty 1_{\{\tau_i\leq nh\}}\beta(\tau_i) K(\rho(X_{\tau_i}^{i},X_{\tau_i^h}^{i})\right] \\
    & =  \mathbb{E}_{x}^{V^{\varepsilon}}\left[\sum_{i=1}^\infty 1_{\{\tau_i\leq nh\}}\beta(\tau_i) \mathbb{E}_{x}^{V^\varepsilon}\left[K(\rho(X_{\tau_i}^{i},X_{\tau_i^h}^{i})|\mathcal{F}_{\tau_i}\right]\right]\\
    & \leq  \mathbb{E}_{x}^{V^{\varepsilon}}\left[\sum_{i=1}^\infty 1_{\{\tau_i\leq nh\}}\beta(\tau_i) \mathbb{E}_{X_{\tau_i}^{i}}\left[\sup_{s\in [0,h]}K(\rho(X_0,X_{s}))\right]\right]\\
    & \leq  \mathbb{E}_{x}^{V^{\varepsilon}}\left[\sum_{i=1}^\infty 1_{\{\tau_i\leq nh\}}\beta(\tau_i)\mathbb{E}_{X_{\tau_i}^i}\left[1_{\{\sup_{s\in [0,h]}\rho(X_0,X_{s})<r\}} h\right]\right]\\
    &\phantom{=} +   \mathbb{E}_{x}^{V^{\varepsilon}}\left[\sum_{i=1}^\infty 1_{\{\tau_i\leq nh\}}\beta(\tau_i) \mathbb{E}_{X_{\tau_i}^i}\left[1_{\{\sup_{s\in [0,h]}\rho(X_0,X_{s})\geq r\}} \Vert K\Vert\right]\right]\\
    & \leq \mathbb{E}_{x}^{V^{\varepsilon}}\left[\sum_{i=1}^\infty 1_{\{\tau_i\leq nh\}}\beta(\tau_i)  \right](h+\gamma(h)\Vert K\Vert).
\end{align*}
Using a similar argument and the fact $\beta(\tau_i^h)\leq \beta(\tau_i)$, we also get
\begin{multline*}
    \mathbb{E}_{x}^{V^{\varepsilon}}\left[\sum_{i=1}^\infty 1_{\{\tau_i\leq nh\}}\beta(\tau_i^h)\left|c(X_{\tau_i^h}^i,X_{\tau_i^h}^{i+1})-c(X_{\tau_{i}^h}^i,X_{\tau_i}^{i+1})  \right| \right]\\
    \leq \mathbb{E}_{x}^{V^{\varepsilon}}\left[\sum_{i=1}^\infty 1_{\{\tau_i\leq nh\}}\beta(\tau_i)  \right](h+\gamma(h)\Vert K\Vert).
\end{multline*}
Next, using~\eqref{eq:beta_mult}, we get $\beta(\tau_i)-\beta(\tau_i^h)\leq \beta(\tau_i)-\beta(\tau_i+h)\leq \beta(\tau_i)(1-\beta(h))$. Thus, we obtain
\begin{multline*}
    \mathbb{E}_{x}^{V^{\varepsilon}}\left[\sum_{i=1}^\infty 1_{\{\tau_i\leq nh\}}(\beta(\tau_i)-\beta(\tau_i^h))c(X_{\tau_{i}^h}^i,X_{\tau_i}^{i+1})   \right]\\
    \leq \left[\sum_{i=1}^\infty 1_{\{\tau_i\leq nh\}} \beta(\tau_i)\right]\Vert c\Vert (1-\beta(h)).
\end{multline*}
Consequently, recalling~\eqref{eq:th:approx_dyadic_disc:1.5}, we have
\begin{multline*}
    \left|\mathbb{E}_{x}^{V^{\varepsilon}}\left[\sum_{i=1}^\infty 1_{\{\tau_i\leq nh\}}\beta(\tau_i)c(X_{\tau_i}^i,X_{\tau_i}^{i+1}) \right]- \mathbb{E}_{x}^{V^{\varepsilon}}\left[\sum_{i=1}^\infty 1_{\{\tau_i^h\leq nh\}}\beta(\tau_i^h)c(X_{\tau_{i}^h}^i,X_{\tau_i^h}^{i+1})  \right]\right|\\
    \leq \mathbb{E}_{x}^{V^{\varepsilon}}\left[\sum_{i=1}^\infty 1_{\{\tau_i\leq nh\}}\beta(\tau_i)  \right](2h+2\gamma(h)\Vert K\Vert+\Vert c\Vert (1-\beta(h))).
\end{multline*}
Combining this with~\eqref{eq:th:approx_dyadic_disc:1} and~\eqref{eq:th:approx_dyadic_disc:1.2}, we obtain
\begin{align}\label{eq:th:approx_dyadic_disc:2}
    0&\leq \inf_{V\in \mathbb{V}_h} J^d(x,V)-\inf_{V\in \mathbb{V}} J^d(x,V) \nonumber\\
    & \leq  \limsup_{n\to\infty} \frac{\mathbb{E}^{V^{\varepsilon}}_x\left[\sum_{i=1}^{\infty}   1_{\{\tau_i\leq nh\}}\beta(\tau_i)\right]}{\int_0^{nh} \beta(s)ds}\left(2\Vert g\Vert h+2h+2\gamma(h)\Vert K\Vert+\Vert c\Vert (1-\beta(h))\right) +\varepsilon.
\end{align}
Next, note that for the $\varepsilon$-optimal strategy $V^\varepsilon$ we have

\[
J(x,V^{\varepsilon})\leq \inf\limits_{V\in \mathbb{V}}J(x,V)+\varepsilon\leq \Vert g\Vert +\varepsilon.
\]
Thus, as in~\eqref{eq:th:approx_disc:2}, we get 
\[
\limsup_{n\to\infty} \frac{\mathbb{E}^{V^{\varepsilon}}_x\left[\sum_{i=1}^{\infty}   1_{\{\tau_i\leq nh\}}\beta(\tau_i)\right]}{\int_0^{nh} \beta(s)ds}\leq \frac{2\Vert g\Vert+\varepsilon}{c_0}.
\]
Consequently, letting $h\downarrow 0$ in~\eqref{eq:th:approx_dyadic_disc:2} and using the fact that $\varepsilon>0$ was arbitrary, we conclude the proof.
\end{proof}
\begin{remark}
    To ensure that the sequence defined in~\eqref{eq:th:approx_dyadic_disc:strategy} is a proper admissible strategy, we used the assumption that we are allowed to shift the controlled process to any point of the state space $E$. Without this assumption, one might proceed as in~\cite{JelPitSte2019b}, where a similar result on the approximation of continuous time impulse control problem with the help of its discrete (dyadic) version is established in the undiscounted risk-sensitive framework; see Theorem 4.4 and the following discussion therein for details. The present paper obtains similar approximation results without technical difficulties related to solving the continuous time Bellman equations.
\end{remark}

In the next theorem, we show that the discrete approximation holds also for the undiscounted problem. The argument is similar to the one used in Theorem~\ref{th:approx_dyadic_disc} and is omitted for brevity.

\begin{theorem}\label{th:approx_dyadic_undisc}
Let $J$ be given by~\eqref{eq:funct_time_cons}. Then, we get 
    \[\lim_{h\to 0}\inf_{V\in \mathbb{V}_h} J(x,V)=\inf_{V\in \mathbb{V}} J(x,V).
    \]
\end{theorem}

Now, we are ready to prove the main result of this paper, stating that the optimal values of discounted and undiscounted problems defined in~\eqref{eq:goal_cont} are equal to each other.

\begin{theorem}\label{th:main}
Let $J^d$ and $J$ be given by~\eqref{eq:funct_time_inc} and~\eqref{eq:funct_time_cons}, respectively. Then, we get 
    \[\inf_{V\in \mathbb{V}} J^d(x,V)=\inf_{V\in \mathbb{V}} J(x,V).
    \]
\end{theorem}
\begin{proof}
    Using successively Theorem~\ref{th:approx_dyadic_disc}, Theorem~\ref{th:approx_disc}, Theorem~\ref{th:main_discrete}, Theorem~\ref{th:approx_undisc}, and Theorem~\ref{th:approx_dyadic_undisc}, we get
    \begin{multline*}
        \inf_{V\in \mathbb{V}} J^d(x,V) = \lim_{h\to 0}\inf_{V\in \mathbb{V}_h} J^d(x,V)  =\lim_{h\to 0}\inf_{V\in \mathbb{V}_h} J^d_h(x,V) \\
        = \lim_{h\to 0}\inf_{V\in \mathbb{V}_h} J_h(x,V) = \lim_{h\to 0}\inf_{V\in \mathbb{V}_h} J(x,V) = \inf_{V\in \mathbb{V}} J(x,V),
    \end{multline*}
    which concludes the proof.
\end{proof}

Now, we show that the strategy, which is optimal for the undiscounted discrete time problem, is also nearly optimal for both discounted and undiscounted continuous time problems.

\begin{corollary}\label{cor:nearly_optimal}
Let $\varepsilon>0$ and $x\in E$. Then, there exists $h=h(\varepsilon)>0$ such that the strategy $\hat{V}_h\in \mathbb{V}_h$ given by~\eqref{eq:opt_str_disc_cons} is $\varepsilon$-optimal for~\eqref{eq:funct_time_inc} and~\eqref{eq:funct_time_cons}.
\end{corollary}
\begin{proof}
    First, using Theorem~\ref{th:approx_dyadic_undisc}, we may find $h_1=h_1(\varepsilon)>0$ such that for $h<h_1$ we get $\inf_{V\in \mathbb{V}_h}J(x,V)\leq \inf_{V\in \mathbb{V}}J(x,V)+\varepsilon/2 $. Next, using Theorem~\ref{th:approx_undisc}, we may find $h_2=h_2(\varepsilon_1)<h_1$, such that for $h<h_2$, we get $\inf_{V\in \mathbb{V}_h}J_h(x,V_h)\leq \inf_{V\in \mathbb{V}_h}J(x,V)+\varepsilon/2 $. Finally, recalling that, by Corollary~\ref{cor:opt_str_disc_undisc}, the strategy  $\hat{V}_h\in \mathbb{V}_h$ is optimal for $\inf_{V\in \mathbb{V}_h}J_h(x,V)$, $x\in E$, we obtain
    \[
J_h(x,\hat{V}_h)=\inf_{V\in \mathbb{V}_h}J_h(x,V)\leq \inf_{V\in \mathbb{V}}J(x,V)+\varepsilon,
    \] 
    which shows that $\hat{V}_h$ is $\varepsilon$-optimal for~\eqref{eq:funct_time_cons}.

    Second, using Theorem~\ref{th:main_discrete}, Theorem~\ref{th:approx_dyadic_disc}, and a similar argument, we may find $h_3=h_3(\varepsilon)<h_2$ such that for $h<h_3$ we get
    \[
\inf_{V\in \mathbb{V}_h}J_h^d(x,V)\leq \inf_{V\in \mathbb{V}}J^d(x,V)+\varepsilon.
    \] 
    Recalling that by Theorem~\ref{th:main_discrete}, the strategy $\hat{V}_h$ is also optimal for $\inf_{V\in \mathbb{V}_h}J_h^d(x,V)$, $x\in E$, we conclude the proof.
\end{proof}

\bibliographystyle{agsm}
\bibliography{RSC_bibliografia}
\end{document}